\documentclass[hidelinks,onefignum,onetabnum]{siamart220329}

\usepackage{amssymb}
\usepackage{amsfonts}
\usepackage{dsfont}
\usepackage{mathrsfs}
\usepackage{graphicx}

\DeclareMathOperator*{\esssup}{ess\,sup}

\usepackage{enumerate}
\usepackage{hyperref}
\usepackage{indentfirst}
\usepackage{color}
\usepackage{multicol,colordvi,epsfig,float}
\usepackage{pgfplots}
\pgfplotsset{width=9cm,compat=1.15}

\newcommand{\dist}[2]{\langle \ {#1} \ , \ {#2} \ \rangle}

\newcommand{\norm}[2]{| {#1} | _{#2}}
\newcommand{\norma}[2]{\| {#1} \|_{#2}}
\newcommand{\prodl}[2]{\left (\ {#1} \ , \ {#2} \ \right )}

\newcommand{\D}[1]{\displaystyle{#1}}

\usepackage{verbatim}

\ifpdf
  \DeclareGraphicsExtensions{.eps,.pdf,.png,.jpg}
\else
  \DeclareGraphicsExtensions{.eps}
\fi

\newsiamremark{remark}{\bf Remark}
\newsiamremark{hypothesis}{Hypothesis}
\crefname{hypothesis}{Hypothesis}{Hypotheses}
\newsiamthm{claim}{Claim}

\headers{Optimal control for chemotaxis-consumption models}{F. Guillén-González, A. L. Corrêa Vianna Filho}

\title{An optimal control problem subject to strong solutions of chemotaxis-consumption models \thanks{
This work was partially funded by Grant PGC2018-098308-B-I00 (MCI/AEI/FEDER, UE). FGG has also been financed in part by the Grant US-1381261 (US/JUNTA/FEDER, UE) and Grant P20-01120 (PAIDI/JUNTA/FEDER,UE).}}   

\author{Francisco Guillén-González\thanks{Departament of Partial Differential Equations and Numerical Analysis, Universidad de Sevilla, Seville, Spain. (\email{guillen@us.es}, \email{acorreaviannafilho@us.es}) .}
\and André Luiz Corrêa Vianna Filho\footnotemark[2]}

\usepackage{amsopn}

\begin{document}

\maketitle

\begin{abstract}
  We consider a bilinear optimal control problem associated to the following chemotaxis-consumption model in a bounded domain $\Omega \subset \mathbb{R}^3$ during a time interval $(0,T)$: $$\partial_t u - \Delta u  = - \nabla \cdot (u \nabla v), \quad \partial_t v - \Delta v  = - u^s v + f v 1_{\Omega_c},$$ with $s \geq 1$, endowed with isolated boundary conditions and initial conditions for $(u,v)$,  $u$ being the cell density, $v$ the chemical concentration and $f$ the bilinear control acting in a subdomain $\Omega_c \subset \Omega$. The existence of weak solutions $(u,v)$ to this model given $f \in L^q((0,T) \times \Omega)$, for some $q > 5/2$, has been proved in \cite{guillen2022optimal}. In this paper, we study a related optimal control problem in the strong solution setting. First, imposing the regularity criterion $u ^s \in L^q((0,T) \times \Omega)$ ($q > 5/2$) for a given weak solution, we prove existence and uniqueness of global-in-time strong solutions. Then,  the existence of a global optimal solution can be deduced. Finally, using a Lagrange multipliers theorem, we establish first order optimality conditions for any local optimal solution, proving  existence, uniqueness and regularity of the associated Lagrange multipliers.
\end{abstract}

\begin{keywords}
  chemotaxis, consumption, optimal control, regularity criterion, Lagrange multiplier.
\end{keywords}

\begin{MSCcodes}
  49J20, 49K20, 49N60, 35Q92, 92C17.
\end{MSCcodes}

\section{Introduction}

  Chemotaxis is the directed movement of cells induced by the gradient of a chemical substance. The introduction of one of the first mathematical models for chemotaxis is attributed Keller and Segel in two works from 1970 and 1971 \cite{keller1970initiation,keller1971model} which are also regarded by some authors as a development of the work of Patlak \cite{patlak1953random}. Since then, the research on this topic gave rise to different related models, such as models with chemoattraction or chemorepulsion, combined with production or consumption of the chemical substance, with the presence of a logistic growth term, models for angiogenesis, haptotaxis and so on, covering a wide variety of applications of practical interest. From the mathematical point of view, the aforementioned models possess interesting and challenging features that attracted the attention of many authors along the years and make these models still relevant nowadays \cite{bellomo2015toward,horstmann20031970,horstmann20041970}.
  
  Let $\Omega$ be a bounded domain of $\mathbb{R}^3$, denoting by $\Gamma$ its boundary, and define $Q: = (0,T) \times \Omega$, for a fixed given $T > 0$. Let $u=u(t,x)$ and $v=v(t,x)$ be the density of cell population and the concentration of chemical substance, respectively, defined on $(t,x) \in Q$. In the present work we are going to study an optimal control problem related to the following chemotaxis-consumption model:
  \begin{equation} \label{problema_P}
    \left\{\begin{array}{l}
      \partial_t u - \Delta u  = - \nabla \cdot (u \nabla v), \quad
      \partial_t v - \Delta v  = - u^s v, \\
      \partial_{\boldsymbol{n}} u |_{\Gamma}  =  \partial_{\boldsymbol{n}} v |_{\Gamma} = 0, \quad
      u(0) = u^0, \quad v(0) = v^0,
    \end{array}\right.
  \end{equation}
  where $s \geq 1$ is a fixed real number, $\partial_{\boldsymbol{n}} u |_{\Gamma}$ denotes the normal derivative of $u$ on the boundary and $u^0 ,v^0 \geq 0$ are the initial data.
  
  Then, starting from this uncontrolled problem \eqref{problema_P}, we consider a  bilinear control acting on the chemical equation, arriving at the controlled problem:
  \begin{equation} \label{problema_P_controlado}
    \left\{\begin{array}{l}
      \partial_t u - \Delta u  = - \nabla \cdot (u \nabla v), \quad
      \partial_t v - \Delta v  = - u^s v + f v 1_{\Omega_c}, \\
      \partial_{\boldsymbol{n}} u |_{\Gamma}  =  \partial_{\boldsymbol{n}} v |_{\Gamma} = 0, \quad
      u(0)  = u^0, \quad v(0) = v^0,
    \end{array}\right.
  \end{equation}
  where $f$ is the control and $1_{\Omega_c}$ is the characteristic function of the control domain $\Omega_c \subset \Omega$. 
  
  Next we recall some developments in the theory of the chemotaxis model \eqref{problema_P}. In \cite{tao2012eventual} some results can be found about existence of global weak solutions for \eqref{problema_P} with $s = 1$ in $3$D smooth and convex domains. These solutions become smooth after a sufficient large time and their asymptotic behavior is analyzed. For the corresponding parabolic-elliptic simplification of \eqref{problema_P}, also with $s = 1$, the existence and uniqueness of a global classical solution and its asymptotic behavior is studied  in \cite{tao2019global}. Still for $s = 1$, there are also studies coupling  \eqref{problema_P} with fluids, namely, with the (Navier)-Stokes equations \cite{winkler2012global,winkler2014stabilization,jiang2015global,winkler2016global,winkler2017far}.

  Recently, in \cite{ViannaGuillen2023uniform}, the existence results of \cite{tao2012eventual} were extended to more general $3$D and $2$D domains that are neither necessarily smooth nor convex. Moreover, motivated by \cite{guillen2021chemorepulsion,guillen2020study}, where the production term $u$ of a chemorepulsion-production model studied in \cite{cieslak2006global} was generalized to $u^s$, for $s \in (1,2]$, the authors in \cite{ViannaGuillen2023uniform} generalize the consumption term $-uv$ to $-u^sv$ varying the exponent $s \geq 1$.
  
  An interesting and challenging feature of chemotaxis models, is that the $L^{\infty}$-norm of the cell density $u$ may blow up in finite time. Concerning the model \eqref{problema_P} with $s = 1$, existence and uniqueness of a global classical solution that is uniformly bounded up to infinity time is proved in \cite{tao2012eventual} for $2$D smooth and convex domains. For $s \geq 1$ and more general nonconvex $2$D domains, existence and uniqueness of a global strong solution that is uniformly bounded up to infinity time is proved in \cite{ViannaGuillen2023uniform}. 
  To the best of our knowledge, whether or not there exist weak solutions of \eqref{problema_P} in $3$D domains such that $u$ blows up in finite time remains as an open question. Studying conditions that could lead to no-blow-up results for \eqref{problema_P}, with $s = 1$, some authors were able to prove existence and uniqueness of a global classical solution  uniformly bounded up to infinity time in bounded smooth $N$-dimensional domains under the assumption of adequate constraints relating the chemotaxis coefficient and $\norma{v^0}{L^{\infty}(\Omega)}$ \cite{tao2011boundedness,baghaei2017boundedness,fuest2019analysis,frassu2021boundedness}.
  
  Another kind of result consists in giving sufficient conditions about the boundedness of local in time classical solutions of chemotaxis models which avoid blow up at finite time, see for instance  \cite[Lemma 3.2]{bellomo2015toward}. In this direction, the present paper gives a regularity criterion in Theorem \ref{teo_dependencia_em_relacao_ao_controle} below. In particular, Theorem \ref{teo_dependencia_em_relacao_ao_controle} guarantees that if $(u,v)$ is a weak solution of the uncontrolled problem \eqref{problema_P} (taking $f = 0$) and $u$ satisfies the additional regularity $u^s \in L^q(0,T;L^q(\Omega))$, for some $q > 5/2$, then $u$ does not blow up at finite time.
  
  Now, we review some works dedicated to the optimal control problem constrained to chemotaxis models. In $1$D domains for Keller-Segel model or in $2$D domains for other chemotaxis models, where one has the existence and uniqueness of strong solution to the controlled model, it is usual to study the existence of global optimal solution and to derive an optimality system, establishing existence and regularity of Lagrange multipliers for any local optimum. Some references are, \cite{ryu2001optimal} for a Keller-Segel model with distributed control;  \cite{guillen2020optimal,guillen2020bi} for a chemorepulsion-production model; \cite{silva2022bilinear} for a Keller-Segel logistic model;  \cite{yuan2022optimal} for a chemotaxis model with indirect consumption; and \cite{tang2022optimal} for a chemotaxis-haptotaxis model. For $3$D domains this analysis is more complex, mainly because, despite in many cases one has results of existence of weak solutions, there is not uniqueness. To overcome this difficulty, a regularity criterion is introduced, which is a mild additional regularity hypothesis on a weak solution, sufficient to conclude that this weak solution is actually the unique strong solution. For a motivated introduction of this kind of adaptation we refer the reader to \cite{casas1998optimal}, where  an optimal control problem related to the Navier-Stokes equations in $3$D domains is studied. For chemotaxis related works in $3$D domains, we cite \cite{guillen2020regularity} for a chemorepulsion-production model, and \cite{lopez2021optimal} for a chemotaxis-Navier-Stokes-consumption model.

  Although optimal control problems related to chemotaxis models have been the focus of a series of recent works cited above, we still have a relative low number of studies on optimal control problems related to the chemotaxis-consumption model \eqref{problema_P}. Indeed, as far as we know, we can cite \cite{guillen2022optimal}, where  an optimal control problem related to the chemotaxis-consumption model \eqref{problema_P} is studied  in a weak solution setting, and \cite{lopez2021optimal}, where a regularity criterion is proved and its application to a control problem subject to a chemotaxis-Navier-Stokes model is carried out.
  
  Accounting for the exposed so far, the objective of the present work is to study an optimal control problem related to \eqref{problema_P_controlado} in the strong solution setting, establishing the existence of global optimal solution and first order optimality conditions for any local optimal solution, proving the existence, uniqueness and regularity of the associated  Lagrange multipliers. In order to achieve it, we begin by proving a regularity criterion for the controlled problem \eqref{problema_P_controlado}.

  
\subsection{Main contributions of the paper}
  
  Throughout this work we assume that
  \begin{equation} \label{assumptions_set}
    \begin{array}{c}
      \Omega \subset \mathbb{R}^3 \mbox{ is a bounded domain with boundary } \Gamma \mbox{ of class } C^{2,1}, \\
      \Omega_c \subset \Omega \mbox{ is a subdomain with boundary } \Gamma_c \mbox{ locally Lipschitz}, \\
      s \geq 1 \mbox{ and } q > 5/2 \mbox{ are fixed real numbers}.
    \end{array}
  \end{equation}
  
  The first main contribution is to give a regularity criterion that, under a mild additional regularity hypothesis over the $u$-component of a weak solution (see Definition \ref{defi_weak_solution} below) of the controlled problem \eqref{problema_P_controlado}, allows us to conclude that it is actually the unique strong solution of \eqref{problema_P_controlado} (see Definition \ref{defi_strong_solution} below). In this result it is also established the continuous dependence in the strong regularity (see relation \eqref{estimativa_solucao_forte_em_relacao_a_f} below), which is essential to prove the existence of global optimal solution.
 
  Let $X_p$ be the Banach space
  \begin{equation*}
    X_p = \{ v \in L^p(0,T;W^{2,p}(\Omega)) \ : \ \partial_t v \in L^p(Q) \}.
  \end{equation*}

  \begin{remark}
    The space $X_p$ is continuously embedded in $C([0, T];W^{2-2/p,p}(\Omega))$ (see \cite[Theorem III.4.10.2]{amann1995linear}).
    \hfill $\square$
  \end{remark}

  In the sequel, we introduce the concepts of weak and strong solution of \eqref{problema_P_controlado}. We begin by the definition of weak solution, whose existence is proved in \cite{guillen2022optimal}, based in results of \cite{ViannaGuillen2023uniform} for the uncontrolled problem \eqref{problema_P} ($f=0$).
  
  \begin{definition}{\bf (Weak solution of \eqref{problema_P_controlado})} \label{defi_weak_solution}
    Let $s \geq 1$, $q > 5/2$. Let $f \in L^q(Q)$ and $(u^0, v^0) \in L^p(\Omega) \times W^{2-2/q,q}(\Omega)$, with $p = 1 + \varepsilon$, for some $\varepsilon > 0$, if $s = 1$, and $p = s$, if $s > 1$, be non-negative functions. A pair $(u,v)$ is called a weak solution of \eqref{problema_P_controlado} if $u(t,x),v(t,x) \geq 0$ $a.e.$ $(t,x) \in Q$, $(u,v)$ has the regularity,
 
 \noindent for $s \geq 1$,
    \begin{equation*}
      \begin{array}{c}
        u^s \in L^{\infty}(0,T;L^1(\Omega)) \cap L^{5/3}(Q), \\
        v \in L^{\infty}(Q) \cap L^{\infty}(0,T;H^1(\Omega)) \cap L^4(0,T;W^{1,4}(\Omega)) \cap L^2(0,T;H^2(\Omega)),
       \ \partial_t v \in L^{5/3}(Q),
      \end{array}
    \end{equation*}
    for $s \in [1,2)$,
    \begin{equation*}
      u \in L^{5s/(3 + s)}(0,T;W^{1,5s/(3 + s)}(\Omega)), \partial_t u \in L^{5s/(3 + s)}(0,T;(W^{1,5s/(4s-3)}(\Omega))'),
    \end{equation*}
    for $s \geq 2$,
    \begin{equation*}
      u \in L^2(0,T;H^1(\Omega)), \ \partial_t u \in L^2(0,T;(H^1(\Omega))'),
    \end{equation*}
    and satisfies the initial conditions for $(u,v)$, the $u$-equation of \eqref{problema_P_controlado} and the boundary condition of $u$ in the variational sense
    \begin{equation} \label{eq_u_sol_fraca}
      \int_0^T \int_{\Omega} \partial_t u \ \varphi \ dx \ dt + \int_0^T \int_{\Omega} \nabla u \cdot \nabla \varphi \ dx \ dt = \int_0^T \int_{\Omega} u \nabla v \cdot \nabla \varphi \ dx \ dt,
    \end{equation}
    for all $\varphi \in L^{5s/(4s-3)}(0,T;W^{1,5s/(4s-3)}(\Omega))$, if $s \in [1,2)$, and $\varphi \in L^2(0,T;H^1(\Omega))$, if $s \geq 2$, the $v$-equation $a.e.$ $(t,x) \in Q$ (in fact, the $v$-equation is satisfied in $L^{5/3}(Q)$) and, since $\Delta v \in L^2(Q)$, the boundary condition of $v$ in the sense of $H^{-1/2}(\Gamma)$.
    \hfill $\square$
  \end{definition}
  
  \begin{definition}{\bf (Strong solution of \eqref{problema_P_controlado})} \label{defi_strong_solution}
    Let $s \geq 1$, $q > 5/2$. Let $f \in L^q(Q)$ and $u^0, v^0 \in W^{2-2/q,q}(\Omega)$ non-negative functions. A pair $(u,v)$ is called a strong solution of \eqref{problema_P_controlado} if $u(t,x),v(t,x) \geq 0$ $a.e.$ $(t,x) \in Q$, with regularity $(u,v) \in X_q \times X_q$ and satisfying the initial and boundary conditions of \eqref{problema_P_controlado}, the $u$-equation and the $v$-equation of \eqref{problema_P_controlado} $a.e.$ $(t,x) \in Q$. Moreover, since $\Delta u, \Delta v \in L^q(Q)$, $u$ and $v$ satisfy the boundary conditions in the sense of $W^{1-1/q,q}(\Gamma)$ (see \cite[Theorem 1.6]{girault2012finite}).
    \hfill $\square$
  \end{definition}
  
  \begin{remark}
    Since $q > 5/2$, if $(u,v)$ is a strong solution of \eqref{problema_P_controlado} with control $f \in L^q(Q)$ then, in particular, $u,v \in L^{\infty}(Q)$. Then, through a comparison argument we can prove that, for each fixed $f \in L^q(Q)$, the strong solution of \eqref{problema_P_controlado} is unique. We refer the reader to the proof of uniqueness in $2$D domains made in \cite{ViannaGuillen2023uniform} that, in view of the regularity of the strong solution, can be adapted to $3$D domains.
    \hfill $\square$
  \end{remark}
  
  Now we are in position to state the first main result of this paper.

  \begin{theorem}[\bf Regularity criterion] \label{teo_dependencia_em_relacao_ao_controle}
    Assume \eqref{assumptions_set}. Let $(u,v)$ be a weak solution of problem \eqref{problema_P_controlado} with
     $f \in L^q(Q)$. If, additionally, we suppose that $u^s \in L^q(Q)$, then $(u,v) \in X_q \times X_q$ is the unique strong solution of problem \eqref{problema_P_controlado}. Moreover, there is $\mathcal{K} = \mathcal{K}(\norma{u^s}{L^q(Q)},\norma{f}{L^q(Q)}) > 0$, where $\mathcal{K}(\cdot,\cdot)$ is a  continuous and increasing function with respect to each entry, $\norma{u^s}{L^q(Q)}$ and $\norma{f}{L^q}$, such that
    \begin{equation} \label{estimativa_solucao_forte_em_relacao_a_f}
      \norma{(u,v)}{X_q \times X_q} \leq \mathcal{K}(\norma{u^s}{L^q(Q)},\norma{f}{L^q(Q)}).
    \end{equation}
  \end{theorem}
  
  \begin{remark}
    Following the proof of Theorem \ref{teo_dependencia_em_relacao_ao_controle} we observe that the power $5/2$ is critical in the sense that Theorem~\ref{teo_dependencia_em_relacao_ao_controle} is proved for any $q > 5/2$ and, at least using the techniques employed in this proof, it is not possible to reach the same conclusion if $q \leq 5/2$. We also note that the hypothesis $f \in L^q(Q)$ with $q > 5/2$ is essential in the proof of existence of weak solutions of \eqref{problema_P_controlado} given in \cite[Lemma 3.1]{guillen2022optimal}. Moreover, since $q > 5/2$ then $X_q \hookrightarrow L^{\infty}(Q)$ (see Lemma \ref{lema_regularidade_w_in_X_p} below). Then Theorem \ref{teo_dependencia_em_relacao_ao_controle} also gives a regularity hypothesis over a weak solution of the controlled problem \eqref{problema_P_controlado} which avoids blow up at finite time.
    \hfill $\square$
  \end{remark}
  
  The second main contribution of this paper is the existence of optimal solution to the following optimal control problem. Let $\mathcal{F}$ be a closed and convex subset of $L^q(Q)$, for a given $q > 5/2$. Consider the cost functional $J: L^{sq}(Q) \times L^2(Q) \times \mathcal{F} \longrightarrow \mathbb{R}$ given by
  \begin{equation} \label{funcional_J}
    \begin{array}{l}
      J(u,v,f) : = \dfrac{\gamma_u}{s q} \D \int_0^T{\norma{u(t) - u_d(t)}{L^{sq} }^{sq} \ dt} \\[6pt]
      \D + \dfrac{\gamma_v}{2} \int_0^T{\norma{v(t) - v_d(t)}{L^2 }^2 \ dt} + \dfrac{\gamma_f}{q} \int_0^T{\norma{f(t)}{L^q(\Omega_c)}^q \ dt},
    \end{array}
  \end{equation}
  where $(u_d,v_d) \in L^{sq}(Q) \times L^2(Q)$ represents the desired states and the parameters $\gamma_u, \gamma_v, \gamma_f \geq 0$ measure the costs of the states and control. In addition, we assume 
  \begin{equation} \label{hipotese_sobre_os_custos_gamma}
    \begin{array}{c}
      \gamma_u > 0 \ \mbox{ and} \\
      \gamma_f > 0 \ \mbox{ or } \ \mathcal{F} \mbox{ is bounded in } L^q(Q).
    \end{array}
  \end{equation}
  
  We are going to minimize $J(u,v,f)$ subject to the admissible set of the triples $(u,v,f)$ satisfying the controlled problem \eqref{problema_P_controlado} in the strong setting, that is 
  \begin{equation*}
    \begin{array}{rcl}
      S_{ad}  =  \{ (u,v,f) \in X_q \times X_q \times \mathcal{F} \ | \ (u,v) \mbox{ is the strong solution of \eqref{problema_P_controlado} with control } f \}.
    \end{array}
  \end{equation*}
   Then, the following minimization problem is considered:
  \begin{equation} \label{problema_de_minimizacao}
    min \ J(u,v,f) \mbox{ subject to } (u,v,f) \in S_{ad}.
  \end{equation}
  Since given $f\in \mathcal{F}$ one can not assure in general the existence of a strong solution $(u,v)$ associated to $f$,  we fix the hypothesis
  \begin{equation} \label{hipotese_conjunto_admissivel_nao_vazio}
    S_{ad} \neq \emptyset.
  \end{equation}
  \begin{remark}
   Analogously to \cite{guillen2020regularity} and \cite{lopez2021optimal}, if $\Omega_c = \Omega$, that is, if the control acts in the whole domain, then  \eqref{hipotese_conjunto_admissivel_nao_vazio} holds. In addition, if we assume that $\Omega$ is a $2$D domain and $0 \in \mathcal{F}$ then \eqref{hipotese_conjunto_admissivel_nao_vazio} also holds. Indeed, from \cite[Theorem 4]{ViannaGuillen2023uniform} we have the existence and uniqueness of weak solution $(u,v)$ with $u \in L^{\infty}(Q)$, of the uncontrolled problem, that is \eqref{problema_P_controlado} with $f = 0$. Since $(u,v)$ and $f=0$ satisfy the hypotheses of Theorem \ref{teo_dependencia_em_relacao_ao_controle}, we conclude that $(u,v) \in X_q \times X_q$ is the strong solution of \eqref{problema_P_controlado} with $f \equiv 0$. In particular, $(u,v,0) \in S_{ad}$ and hence $S_{ad} \neq \emptyset$.
    \hfill $\square$
  \end{remark}
  \begin{theorem}[\bf Existence of optimal control] \label{teo_existencia_controle_otimo}
    Recall the hypotheses \eqref{assumptions_set} and, moreover, assume $S_{ad} \neq \emptyset$. Then the optimal control problem \eqref{problema_de_minimizacao} has at least one global optimal solution $(\overline{u},\overline{v},\overline{f}) \in S_{ad}$.
  \end{theorem}
   
  The third main contribution of this work is the existence and uniqueness of Lagrange multipliers associated to any local optimal solution of \eqref{problema_de_minimizacao}. Let $q' = q/(q-1)$, the conjugate exponent of $q$, let $ (\overline{u},\overline{v},\overline{f}) \in S_{ad}$ be a local optimal solution of \eqref{problema_de_minimizacao} and consider the following Lagrange multiplier problem for $(\lambda, \eta)$ associated to $(\overline{u},\overline{v},\overline{f})$:
  \begin{equation} \label{problema_adjunto_ao_linearizado}
     \left \{
     \begin{array}{l}
       - \partial_t \lambda - \Delta \lambda - \nabla \overline{v} \cdot \nabla \lambda + s \overline{u}^{s-1} \overline{v} \eta = g_\lambda, \\[6pt]
       - \partial_t \eta - \Delta \eta + \overline{u}^s \eta - \overline{f} \eta \ 1_{\Omega_c} + \nabla \cdot (\overline{u} \nabla \lambda) = g_\eta, \\[6pt]
       \partial_{\boldsymbol{n}} \lambda |_{\Gamma} = \partial_{\boldsymbol{n}} \eta |_{\Gamma} = 0, \ \lambda(T,x) = \eta(T,x) = 0,
     \end{array}
     \right.
   \end{equation}
   where
   \begin{equation} \label{rhs-adjoint}
     g_\lambda= \gamma_u sgn(\overline{u} - u_d) \norm{\overline{u} - u_d}{}^{sq - 1}
     \quad  \hbox{and} \quad
     g_\eta=\gamma_v (\overline{v} - v_d).
   \end{equation}
   \begin{definition}{\bf (Very weak solution of \eqref{problema_adjunto_ao_linearizado})} \label{defi_very_weak_solution}
     Let $s \geq 1$, $q > 5/2$ and $q' = q/(q - 1)$. A pair $(\lambda,\eta) \in L^{q'}(Q) \times L^{q'}(Q)$ is called a very weak solution of \eqref{problema_adjunto_ao_linearizado} if $(\lambda,\eta)$ satisfies \eqref{problema_adjunto_ao_linearizado} in the sense of the dual space of $X_q \times X_q$, that is, the following variational formulation holds for any $U,V \in X_q$ with 
     $\partial_{\boldsymbol{n}} U  |_{\Gamma} = \partial_{\boldsymbol{n}} V  |_{\Gamma}=0$ and $U(0)=V(0)=0$:
     \begin{equation}\label{adjoint-1}
      \begin{array}{l}
      \D  {\int_0^T \!\!\int_{\Omega}} \lambda \Big ( \partial_t U - \Delta U + \nabla \cdot (U \nabla \overline{v}) \Big ) \ dx \ dt + {\int_0^T \!\! \int_{\Omega}} s \overline{u}^{s-1} \overline{v} \eta \ U \ dx \, dt = {\int_0^T \int_{\Omega}} g_\lambda U \ dx \,dt, 
      \end{array}
    \end{equation}
    \vspace{-12pt}
    \begin{equation}\label{adjoint-2}
      \begin{array}{l}
   \D    {\int_0^T\!\! \int_{\Omega}} \eta \Big ( \partial_t V - \Delta V + \overline{u}^s V - \overline{f} V 1_{\Omega_c} \Big ) \ dx \,dt + {\int_0^T \!\!\int_{\Omega}} \lambda \ \nabla \cdot (\overline{u} \nabla V) \ dx \,dt = \D{\int_0^T \!\!\int_{\Omega}} g_\eta \ V \ dx \, dt.
     \end{array}
    \end{equation}
     \hfill $\square$
   \end{definition}
   
   \begin{theorem}[\bf Existence of Lagrange multipliers] \label{teo_existencia_multiplicador_de_Lagrange_aplicado}
    Assume \eqref{assumptions_set} and let $(\overline{u},\overline{v},\overline{f}) \in S_{ad}$ be a local optimal solution of \eqref{problema_de_minimizacao}. Then there exists a unique Lagrange multiplier $(\lambda,\eta) \in L^{q'}(Q) \times L^{q'}(Q)$ which is a very weak solution of the optimality system \eqref{problema_adjunto_ao_linearizado} and the following optimality condition holds:
    \begin{equation} \label{condicao_otimalidade_fraco_F}
      \D{\int_0^T \int_{\Omega_c}} (\gamma_f sgn(\overline{f}) \norm{\overline{f}}{}^{q-1} + \overline{v}\, \eta)(f - \overline{f}) \ dx \ dt \geq 0, \quad \forall f \in \mathcal{F}.
    \end{equation}
  \end{theorem}
  \begin{remark}
    If $\gamma_f > 0$ and there is no convex constraint on the control, that is $\mathcal{F} = L^q(Q)$, then \eqref{condicao_otimalidade_fraco_F} is equivalent to
    $  \gamma_f sgn(\overline{f}) \norm{\overline{f}}{}^{q-1} + \overline{v} \,\eta = 0$.
    Since $\overline{v} \geq 0$, we conclude the following explicit expression for the control $\overline{f} = - sgn(\eta) \left ( \dfrac{1}{\gamma_f} \overline{v}\, \norm{\eta}{} \right )^{1/(q-1)}$.
    \hfill $\square$
  \end{remark}

  The key to establish the existence of a Lagrange multiplier is to prove the existence of solution to the linearized problem given in \eqref{problema_linearizado} below. To help with this proof, in the Appendix, we provide a result of existence of solution to an adequate general parabolic linear system. This result is also useful in the study of the regularity of the Lagrange multiplier $(\lambda, \eta)$ provided by Theorem \ref{teo_existencia_multiplicador_de_Lagrange_aplicado} depending on the $L^p$ regularity of the RHS term $g_{\lambda}$ given in \eqref{rhs-adjoint}. 
      
  \begin{theorem} \label{teo_regularidade_adicional_multiplicadores_de_Lagrange}
    Assume \eqref{assumptions_set} and let $(\overline{u},\overline{v},\overline{f}) \in S_{ad}$ be a local optimal of problem \eqref{problema_de_minimizacao}. It holds:
    \begin{enumerate}
      \item if $g_{\lambda} \in L^p(Q)$, for $p\in [10/9 , 10/7)$, then 
      the Lagrange multiplier $(\lambda,\eta) \in L^2(Q) \times L^2(Q)$ and satisfies  \eqref{problema_adjunto_ao_linearizado} in the very weak sense (as in \eqref{adjoint-1}-\eqref{adjoint-2});
      \item if $g_{\lambda} \in L^p(Q)$, for $p\in [10/7 , 2]$, then 
      the Lagrange multiplier $(\lambda,\eta) \in X_p \times X_p$ and satisfies 
      \eqref{problema_adjunto_ao_linearizado} in the strong sense, that is, $a.e.$ $(t,x)$ in $Q$.
    \end{enumerate}
  \end{theorem}

  \begin{remark}
    Since we consider $v_d \in L^2(Q)$, which implies $g_{\eta} \in L^2(Q)$, the previous analysis for $p > 2$ does not seem to lead to  more relevant conclusions.
    \hfill $\square$
  \end{remark}
  \begin{remark}
    To guarantee that the terms of the functional $J$ given in \eqref{funcional_J} make sense it is enough that $u_d \in L^{\tilde q}(Q)$, with $\tilde{q} \geq sq$, and $v_d \in L^2(Q)$. With this regularity,  $g_{\eta} \in L^2(Q)$ and $g_{\lambda} \in L^p(Q)$, for a power $p = p(s,q,\tilde{q}) = \tilde{q}/(sq - 1)$. Hence the regularity of $g_{\lambda}$ depends on $s \geq 1$, $q > 5/2$ and $\tilde{q} \geq sq$, and is decreasing with respect to $s$, with $p(s,q,\tilde{q}) \to 1$ as $s \to \infty$. For instance if  $\tilde{q} = sq$, we have $p = sq/(sq - 1)$. In this case, since $s \geq 1$ and $q > 5/2$, then $p \in (1,5/3)$. Let us fix $q > 5/2$ close to $5/2$ and vary the values of $s$. Then,  if $s \in [1,10/3q]$  we are in the case $(2)$ of Theorem \ref{teo_regularidade_adicional_multiplicadores_de_Lagrange}, and if $s \in (10/3q,10/q]$  we are in the case $(1)$ of Theorem \ref{teo_regularidade_adicional_multiplicadores_de_Lagrange}. But, if $s > 10/q$ then $p \in (1,10/9)$, hence  Theorem \ref{teo_regularidade_adicional_multiplicadores_de_Lagrange} doesn't give  
 additional regularity for the Lagrange multiplier.
    \hfill $\square$
  \end{remark}
  
  The rest of the paper is organized as follows. In Section \ref{section: preliminary} we introduce some notation and  preliminary results that will be used along this paper. In Section \ref{section: regularity criterion} we establish some previous results and prove Theorem \ref{teo_dependencia_em_relacao_ao_controle}. In Section \ref{section: existencia controle otimo} we prove Theorem \ref{teo_existencia_controle_otimo}. Theorems \ref{teo_existencia_multiplicador_de_Lagrange_aplicado} and \ref{teo_regularidade_adicional_multiplicadores_de_Lagrange} are proved in Section \ref{section: multiplicadores de Lagrange}.

\section{Notation and preliminary results}
  \label{section: preliminary}
  
  Let $X$ and $Y$ be Banach spaces, we say that $X$ is continuously injected in $Y$, and denote it by $X \hookrightarrow Y$, if $X \subset Y$ and, moreover, there is a constant $C > 0$ such that $\norma{\varphi}{Y} \leq C \norma{\varphi}{X}, \ \forall \, \varphi \in X$.
  
  For $p \in [1,\infty]$, we denote by $L^p(\Omega)$, the usual Banach spaces of $p$-integrable Lebesgue-measurable functions, with the norm $\norma{\cdot}{L^p }$. We denote by $p'=p/(p-1)$ the conjugate exponent of $p$. We recall that $L^2(\Omega)$ is a Hilbert space with the inner product
   $ \prodl{f}{g} = \int_{\Omega}{f(x) g(x) \, dx}$. 
  We also denote by $W^{k,p}(\Omega)$, with $k\in \mathbb{N}$,  the usual Sobolev space, equipped with the usual norm $\norma{\cdot}{W^{k,p}(\Omega)}$; for $p=2$, we denote $W^{k,2}(\Omega)$ by $H^{k}(\Omega)$, with norm $\norma{\cdot}{H^{k}(\Omega)}$. 
  
  If $X$ is a Banach space, then $L^p(0,T;X)$ is the Bochner space with the norm
  \begin{equation*}
    \norma{v}{L^p(0,T;X)} = \left ( \int_0^T{\norma{v(t)}{X}^p \ dt} \right )^{1/p}, \quad \norma{v}{L^{\infty}(0,T;X)} = \esssup_{t \in (0,T)} \norma{v(t)}{X}.
  \end{equation*}
  To simplify the notation, from now on, we denote the spaces $L^p(\Omega)$ and $W^{k,p}(\Omega)$ by $L^p$ and $W^{k,p}$, respectively, suppressing the domain $\Omega$. Analogously, the spaces $L^q(0,T,L^p(\Omega))$ and $L^q(0,T,W^{k,p}(\Omega))$ will be denoted by $L^q(0,T,L^p)$ and $L^q(0,T,W^{k,p})$. The spaces $L^p(0,T,L^p(\Omega))$ will keep being denoted by $L^p(Q)$.
  
  If $p = 2$ and $X$ is a Hilbert space then $L^2(X)$ is a Hilbert space with the inner product
  \begin{equation*}
    \prodl{u}{v}_{L^2(X)} = \int_0^T{\prodl{u(t)}{v(t)}_X \ dt}, \quad \forall u,v \in L^2(X),     
  \end{equation*}
  where $\prodl{\cdot}{\cdot}_X$ denotes the inner product of $X$.

  Next we introduce some technical lemmas that will be useful throughout the paper.

  \begin{lemma} \label{lema_interpolacao_norma_10/3_e_10}
    Let $\Omega \subset \mathbb{R}^3$ be a bounded Lipschitz domain. There is a constant $C > 0$ such that
    \begin{equation} \label{interpolacao_norma_10/3}
      \norma{v}{L^{10/3}} \leq C \norma{v}{L^2}^{2/5} \norma{v}{H^1}^{3/5}, \ \forall v \in H^1,
    \end{equation}
    \vspace{-12pt}
    \begin{equation} \label{interpolacao_norma_10}
      \norma{v}{L^{10}} \leq C \norma{v}{H^1}^{4/5} \norma{v}{H^2}^{1/5}, \ \forall v \in H^2.
    \end{equation}
  \end{lemma}
  \begin{proof}[\bf Proof]
    Using the interpolation inequality (see \cite{brezis2011functional})
    \begin{equation} \label{general_interpolation_inequality_in_L^p}
      \norma{f}{L^r} \leq \norma{f}{L^p}^{\theta} \norma{f}{L^q}^{1 - \theta}, \mbox{ where } \frac{1}{r} = \frac{\theta}{p} + \frac{1 - \theta}{q} \mbox{ and } \theta \in [0,1],
    \end{equation}
    with $r = 10/3$, $p = 2$ and $q = 6$ we get $\theta = 2/5$ and then
    \begin{equation*}
      \norma{v}{L^{10/3}} \leq \norma{v}{L^2}^{2/5} \norma{v}{L^6}^{3/5}, \quad \forall v \in L^6.
    \end{equation*}
    Thus, from the Sobolev embedding $H^1 \hookrightarrow L^6$ in $3$D domains (see \cite{brezis2011functional}) we obtain \eqref{interpolacao_norma_10/3}.

    Now we prove \eqref{interpolacao_norma_10}. By applying \eqref{general_interpolation_inequality_in_L^p} with $r = 30/13$, $p = 2$ and $q = 6$ we get $\theta = 4/5$ and the Sobolev embedding $W^{1,30/13} \hookrightarrow L^{10}$ in $3$D domains, one has
    \begin{equation*}
      \norma{v}{L^{10}} \leq  C\, \norma{v}{W^{1,30/13}} \leq C \, \norma{v}{H^1}^{4/5} \norma{v}{W^{1,6}}^{1/5}, \quad \forall v \in W^{1,6}.
    \end{equation*}
    Finally, using the Sobolev embedding $H^2 \hookrightarrow W^{1,6}$ in $3$D domains we obtain \eqref{interpolacao_norma_10}.
  \end{proof}

  \begin{lemma}{\bf (\cite{guillen2020regularity})} \label{lema_interpolacao_L_p_W_m_p}
    Let $\Omega \subset \mathbb{R}^N$ be a bounded Lipschitz domain and let $p_1,q_1,p_2,\tilde{p},\tilde{q} \geq 1$ be such that
    \begin{equation*}
      \frac{1}{\tilde{q}} = \frac{(1 - \theta)}{q_1} + \theta\left (\frac{1}{p_1} - \frac{r}{N} \right), \mbox{ and } \frac{1}{\tilde{p}} = \frac{\theta}{p_2}, \mbox{ with } \theta \in [0,1] \mbox{ and } r > 0,
    \end{equation*}
    then $L^{\infty}(L^{q_1}) \cap L^{p_2}(W^{r,p_1}) \hookrightarrow L^{\tilde{p}}(L^{\tilde{q}}).$
  \end{lemma}

\
  
  \begin{lemma}{\bf (\cite[Section 0.4]{feireisl2009singular})} \label{lema_interpolacao_derivada_fracionaria}
    Let $\Omega \subset \mathbb{R}^N$ be a bounded Lipschitz domain. Then the interpolation inequality
    \begin{equation}
      \norma{w}{W^{\alpha,r} } \leq C \norma{w}{W^{\beta,\tilde{p}} }^{\lambda} \norma{w}{W^{\gamma,\tilde{q}} }^{1 - \lambda}, \quad \forall w \in W^{\beta,\tilde{p}} \cap W^{\gamma,\tilde{q}},
    \end{equation}
    holds for $0 \leq \alpha, \beta, \gamma, \lambda \leq 1$ and $1 < \tilde{p},\tilde{q},r < \infty$ such that $\alpha = \lambda \beta + (1 - \lambda) \gamma$ and $\dfrac{1}{r} = \dfrac{\lambda}{\tilde{p}} + \dfrac{(1 - \lambda)}{\tilde{q}}$.
  \end{lemma}

  \begin{remark}
    The spaces with fractional derivatives $W^{\alpha,r}$, $W^{\beta,\tilde{p}}$ and $W^{\gamma,\tilde{q}}$ are the so called \emph{Sobolev-Slobodeckii} spaces. For more details, we refer the reader to \cite[Section 0.4]{feireisl2009singular} and the references suggested therein.
    \hfill $\square$
  \end{remark}
  
\begin{lemma} \label{lema_regularidade_w_in_X_p}
    Let $\Omega \subset \mathbb{R}^3$ be a bounded Lipschitz domain. It holds:
    \begin{enumerate}
      \item $X_p \hookrightarrow L^{5p/(5-2p)}(Q)$, if $p \in [1,5/2)$;
      \item $X_p \hookrightarrow L^{\infty}(L^q)$, for all $q \in [1,\infty)$, if $p = 5/2$;
      \item $X_p \hookrightarrow L^{\infty}(Q)$ if $p > 5/2$.
    \end{enumerate}
  \end{lemma}
  \begin{proof}[\bf Proof]
    By definition of $X_p$, if $w \in X_p$ then we have $w \in C(W^{2-2/p,p}) \cap L^p(W^{2,p})$. If $p \in [1,5/2)$, this implies $w \in C(W^{2-2/p,p}) \cap L^p(W^{2,p}) \hookrightarrow L^{\infty}(L^{3p/(5-2p)}) \cap L^p(W^{2,p})$. Then Lemma \ref{lema_interpolacao_L_p_W_m_p} yields the desired result for  $p < 5/2$. For $p = 5/2$ we use the continuous injection $W^{2-2/p,p}  \hookrightarrow L^q $, for all $q \in [1,\infty)$, and for $p > 5/2$, the continuous injection $W^{2-2/p,p}  \hookrightarrow L^{\infty} $.
  \end{proof}
  
  \begin{lemma} \label{lema_regularidade_nabla_w_in_X_p}
    Let $\Omega \subset \mathbb{R}^3$ be a bounded Lipschitz domain and $p \in (1,5)$. If $w \in X_p$ then $\nabla w \in L^{5p/(5-p)}(Q)$. Moreover, there is a constant $C > 0$ such that,
    \begin{equation*}
      \norma{\nabla w}{L^{5p/(5-p)}(Q)} \leq C \norma{w}{X_p}, \forall w \in X_p.
    \end{equation*}
  \end{lemma}
  \begin{proof}[\bf Proof]
  
    Case $p \in [2,5)$. Since $w \in X_p$, we have by definition that
    \begin{equation*}
      \nabla w \in L^{\infty}(W^{1-2/p,p}) \cap L^p(W^{1,p}) \hookrightarrow L^{\infty}(L^{3p/(5-p)}) \cap L^p(W^{1,p})
    \end{equation*}
    and using Lemma \ref{lema_interpolacao_L_p_W_m_p} we conclude that $\nabla w \in L^{5p/(5-p)}(Q).$ Case $p \in (1,2)$. From the definition of $X_p$, $w \in L^{\infty}(W^{2-2/p,p}) \cap L^p(W^{2,p}),$ hence $D^{2-2/p} w \in L^{\infty}(L^p) \cap L^p(W^{2/p,p})$ and this implies that
    \begin{equation*}
      D^{2-2/p} w \in L^{\infty}(L^p) \cap L^p(W^{\beta,3 p/(1 + \beta p)}), \mbox{ for any } \beta \in (1,2/p).
    \end{equation*}
    Now, using Lemma \ref{lema_interpolacao_derivada_fracionaria} with
      $\alpha = \frac{2}{p} - 1, \ \beta = \beta, \ \tilde{p} = \frac{3p}{1 + \beta p}, \ \gamma = 0 \mbox{ and } \tilde{q} = p$
    and hence
    \begin{equation*}
      \lambda = \frac{2}{p} - 1 \mbox{ and } r = \frac{3 \beta p^2}{- \beta p^2 + (5 \beta + 2)p - 4},
    \end{equation*}
    we obtain 
    \begin{equation*}
      \norma{D^{2-2/p} w}{W^{2/p - 1,r} }^r \leq C \norma{D^{2-2/p} w}{W^{\beta,3p/(1 + \beta p)} }^{(2/p - 1)r} \norma{D^{2-2/p} w}{L^{p} }^{(2 - 2/p)r}.
    \end{equation*}
    The right hand side of this inequality will be integrable choosing $\beta$ such that $(2/p - 1)r = p$. Therefore, choosing $\beta = \dfrac{10 - 5p}{5p - p^2}$, we conclude that $\nabla w \in L^r(Q)$, with $r = \dfrac{10p - 5p^2}{p^2 - 7p + 10} = \dfrac{5p}{5 - p}.$
  \end{proof}
   
  \begin{lemma}{\bf (Compactness in Bochner spaces, \cite{Simon1986compact})} \label{lema_Simon}
      Let $X,B$ and $Y$ be Banach spaces such that $X \subset B \subset Y$, with compact embedding $X \subset B$ and continuous embedding $B \subset Y$. Let $\boldsymbol{F}$ be a set such that $\boldsymbol{F} \subset \Big \{ f \in L^1(0,T;Y) \ \Big | \ \partial_t f \in L^1(0,T;Y) \Big \}$. We have:
      \begin{enumerate}
        \item if the set $\boldsymbol{F}$ is bounded in $L^q(0,T;B) \cap L^1(0,T;X)$, for $1 < q \leq \infty$, and $\Big \{ \partial_t f, \ \forall f \in \boldsymbol{F} \Big \}$ is bounded in $L^1(0,T;Y)$, then $\boldsymbol{F}$ is relatively compact in $L^p(0,T;B)$, for $1 \leq p < q$;
        \item if $\boldsymbol{F}$ is bounded in $L^{\infty}(0,T;X)$ and $\Big \{ \partial_t f, \ \forall f \in \boldsymbol{F} \Big \}$ is bounded in $L^r(0,T;Y)$ for some $r > 1$, then $\boldsymbol{F}$ is relatively compact in $C([0,T];B)$.
      \end{enumerate}
    \end{lemma}

\

     \begin{lemma}{\bf (\cite[Theorem 10.22]{feireisl2009singular})} \label{lema_regularidade_eq_calor}
    Let $\Omega$ be a bounded domain of $\mathbb{R}^N$ such that $\Gamma$ is of class $C^2$. Let $p \in (1,3)$, $w^0 \in W^{2-2/p,p} $ and $h \in L^p(Q)$. Then the problem
    \begin{equation*}
      \left\{\begin{array}{l}
        \partial_t w - \Delta w = h \mbox{ in } Q, \\
        \partial_{\boldsymbol{n}} w |_{\Gamma}  =  0 \mbox{ on } (0,T) \times \Gamma, \\
        w(0,x)  = w^0(x) \mbox{ in } \Omega,
      \end{array}\right.
    \end{equation*}
    has a unique solution $w \in X_p$. 
    Moreover, there is $C = C(p,T,\Omega)>0$ such that
    \begin{equation} \label{dependencia_continua_de_w_em_h}
        \norma{w}{X_p} 
        \leq C ( \norma{h}{L^p(Q)} + \norma{w^0}{W^{2 - 2/p,p} }).
    \end{equation}
  \end{lemma}

\

  \begin{remark}
    Because of the assumption on $\Omega$ in \eqref{assumptions_set}, we have, in particular, that Lemma \ref{lema_regularidade_eq_calor} is applicable.
    \hfill $\square$
  \end{remark}

  \begin{lemma}  \label{lema_termo_fonte_final}
    Suppose that $\Omega \subset \mathbb{R}^3$ is a bounded domain with boundary $\Gamma$ of class $C^{2,1}$. Then there exist positive constants $C_1, C_2 > 0$ such that
    \begin{align*}
      \int_{\Omega}{\norm{\Delta z}{}^2 \ dx} + \int_{\Omega}{\frac{\norm{\nabla z}{}^2}{z} \Delta z \ dx} & \geq C_1 \Big ( \int_{\Omega}{\norm{D^2 z}{}^2 \ dx} + \int_{\Omega}{\frac{\norm{\nabla z}{}^4}{z^2} \ dx} \Big ) - C_2 \int_{\Omega}{\norm{\nabla z}{}^2 \ dx},
    \end{align*}
    for all $z \in H^2(\Omega)$ such that $\partial_\eta z \Big |_{\Gamma} = 0$ and $z \geq \alpha$, for some $\alpha > 0$.
  \end{lemma}
  \begin{proof}[\bf Proof]
    This result is a consequence of \cite[Lemma 22]{ViannaGuillen2023uniform} and \cite[Appendix A.1]{ViannaGuillen2023uniform}.
  \end{proof}


\section{Regularity Criterion}
  \label{section: regularity criterion}
  
  The main objective of the present section is to demonstrate Theorem~\ref{teo_dependencia_em_relacao_ao_controle}. To do so, we first introduce and prove a series of useful results.
  
  \begin{lemma} \label{lema_maquina_bootstrap}
    Let $(u,v)$ be a weak solution of \eqref{problema_P_controlado} (see Definition \ref{defi_weak_solution}). Suppose, in addition, that
    $$u \in L^p(Q), \hbox{ for some $p > 5/3$, } u^0 \in W^{2 - 2/q,q} \mbox{ and } v \in X_q, \hbox{ for some $q > 5/2$.}$$
    Then, $u \in X_{pq/(p+q)}$ and that there is $C = C(\norma{u}{L^p(Q)},\norma{\nabla u}{L^{5/4}(Q)},\norma{v}{X_q}) > 0$, which is continuous and increasing with respect to each entry, $\norma{u}{L^p(Q)}$, $\norma{\nabla u}{L^{5/4}(Q)}$ and $\norma{v}{X_q}$, such that
    \begin{equation} \label{estimativa_maquina_bootstrap_1}
      \norma{u}{X_{pq/(p+q)}} \leq C(\norma{u}{L^p(Q)},\norma{\nabla u}{L^{5/4}(Q)},\norma{v}{X_q}).
    \end{equation}
    The result is also valid for $p = \infty$ and, in this case, we conclude that $u \in X_q$ with
    \begin{equation} \label{estimativa_maquina_bootstrap_2}
      \norma{u}{X_q} \leq C(\norma{u}{L^{\infty}(Q)},\norma{\nabla u}{L^{5/4}(Q)},\norma{v}{X_q}).
    \end{equation}
  \end{lemma}
  \begin{proof}[\bf Proof]
    The basic idea of the proof is a bootstrapping in the $u$-equation of \eqref{problema_P_controlado} that allows one to arrive at the desired regularity  $u \in X_{pq/(p+q)}$ in a finite number of iterations. We are going to consider the case $p < \infty$ and, with small adaptations, one can follow the same reasoning for $p = \infty$.  The proofs of inequalities \eqref{estimativa_maquina_bootstrap_1} and \eqref{estimativa_maquina_bootstrap_2} come from the fact that all the results used along this proof, such as Lemmas \ref{lema_regularidade_nabla_w_in_X_p} and \ref{lema_regularidade_eq_calor}, for example, give us continuous injections. Indeed, since the number of steps of the procedure of gaining regularity will be finite, one can follow the estimates furnished by Lemmas \ref{lema_regularidade_nabla_w_in_X_p} and \ref{lema_regularidade_eq_calor} each time they are applied and, at the end, conclude \eqref{estimativa_maquina_bootstrap_1} and \eqref{estimativa_maquina_bootstrap_2}. Bearing that in mind, we proceed with the proof of $u \in X_{pq/(p+q)}$, for finite $p$.
    
    Using Lemma \ref{lema_regularidade_nabla_w_in_X_p} for $v$ we conclude that $\nabla v \in L^{5q/(5 - q)}(Q)$. Since $q > 5/2$ we have, in particular, that
    \begin{equation} \label{regularidade_nabla_v_L_5+}
      \mbox{there is } \beta > 1 \mbox{ such that } \nabla v \in L^{5 \beta}(Q).
    \end{equation}
    And since $(u,v)$ is a weak solution of \eqref{problema_P_controlado} we have, in particular, $\nabla u \in L^{5/4}(Q)$. By hypothesis and by the definition of $X_q$ we have $u \in L^p(Q)$, with $p > 5/3$ and $\Delta v \in L^q(Q)$ with $q > 5/2$. Considering these regularities, we have the $u$-equation of \eqref{problema_P_controlado} satisfied in the sense
    \begin{equation} \label{u_equation_strong_sense}
      \partial_t u - \Delta u = - u \Delta v - \nabla u \cdot \nabla v, \ a.e. \ (t,x) \in Q,
    \end{equation}
    where
    \begin{equation} \label{regularidade_u_Delta_v}
      u \Delta v \in L^{pq/(p+q)}(Q), \mbox{ with } \frac{pq}{p+q} > 1,
    \end{equation}
    and
    \begin{equation} \label{regularidade_r_0_termo_ gradiente}
      \nabla u \cdot \nabla v \in L^{r_0}, \mbox{ with } r_0 = \frac{5 \beta}{4 \beta + 1} > 1.
    \end{equation}
    Hence, by applying Lemma \ref{lema_regularidade_eq_calor} to  \eqref{u_equation_strong_sense} we conclude that
     $u \in X_r, \mbox{ with } r = \min\left\{r_0,\frac{pq}{p+q} \right\}  > 1.$
    If $r_0 \geq \frac{pq}{p+q}$ then $r = \frac{pq}{p+q}$ and the proof is finished. Therefore it suffices to assume  $r_0 < \frac{pq}{p+q}$. Since for $u \Delta v$ we already have \eqref{regularidade_u_Delta_v}, we focus on enhancing the regularity of the term $\nabla u \cdot \nabla v$.
    
    In this case, we have $u \in X_{r_0}$. Using Lemma \ref{lema_regularidade_nabla_w_in_X_p} we obtain $\nabla u \in L^{5 r_0/ (5 - r_0)}(Q)$. Considering this regularity and \eqref{regularidade_nabla_v_L_5+}, where $\beta > 1$, and \eqref{regularidade_r_0_termo_ gradiente}, where $r_0 > 1$, the new regularity of $\nabla u \cdot \nabla v$ is $L^{\gamma}(Q)$, with $\gamma = \dfrac{5 \beta}{5 \beta - (\beta - 1)r_0} \ r_0 > \dfrac{5 \beta}{4 \beta + 1} \ r_0.$
    
    Define $\alpha = 5 \beta/(4 \beta + 1)$. Note that $\alpha = r_0 > 1$ and $\gamma > \alpha r_0$. Then, let us define $r_1 = \alpha r_0$. Since $\alpha > 1$, we have $r_1 > r_0 > 1$. Now, if $r_1 < \frac{pq}{p+q}$ then, from Lemma \ref{lema_regularidade_eq_calor}, we have $u \in X_{r_1}$. Proceeding by induction, if we have $\nabla u \cdot \nabla v \in L^{r_{n-1}}$, with $r_{n-1} = \alpha^{n-1} r_0 < \frac{pq}{p+q}$, then we have $u \in X_{r_{n-1}}$ and, by Lemma \ref{lema_regularidade_nabla_w_in_X_p} we obtain $\nabla u \in L^{5 r_{n-1}/ (5 - r_{n-1})}(Q)$. And using again \eqref{regularidade_nabla_v_L_5+}, where $\beta > 1$, and \eqref{regularidade_r_0_termo_ gradiente}, where $r_0 > 1$, the new regularity of $\nabla u \cdot \nabla v$ is $L^{\gamma}(Q)$, with
    \begin{equation*}
      \gamma = \frac{5 \beta}{5 \beta - (\beta - 1)r_{n-1}} \ r_{n-1} > \frac{5 \beta}{4 \beta + 1} \ r_{n-1} = \alpha r_{n-1} = \alpha^n r_0.
    \end{equation*}
    Therefore we can define $r_n = \alpha^n r_0$ and again applying Lemma \ref{lema_regularidade_eq_calor} to \eqref{u_equation_strong_sense} we conclude that
    \begin{equation*}
      u \in X_r, \mbox{ with } r = \min \left\{ \alpha^n r_0,\frac{pq}{p+q} \right\} .
    \end{equation*}
    Since $\alpha > 1$, there is an index $n_0$ such that $\alpha^{n_0} r_0 < \frac{pq}{p+q}$ but $\alpha^{n_0 + 1} r_0 \geq \frac{pq}{p+q}$. Therefore we arrive at  $u \in X_{pq/(p+q)}$. 
  \end{proof}

  \begin{remark}
    As it is observed in the proof of Lemma \ref{lema_maquina_bootstrap}, the estimates \eqref{estimativa_maquina_bootstrap_1} and \eqref{estimativa_maquina_bootstrap_2} come from \eqref{dependencia_continua_de_w_em_h}. Therefore, the constants $C$ appearing in \eqref{estimativa_maquina_bootstrap_1} and \eqref{estimativa_maquina_bootstrap_2} also depend on the initial condition $u^0$. Analogously, we conclude that the constant appearing in \eqref{estimativa_solucao_forte_em_relacao_a_f_nabla_u} of Theorem \ref{teo_criterio_de_regularidade} also depends on $u^0$ and $v^0$. But considering that the initial condition $(u^0,v^0) \in W^{2 - 2/q,q} \times W^{2 - 2/q,q}$ of problem \eqref{problema_P_controlado} is fixed, we suppress this dependence in the notation.
    \hfill $\square$
  \end{remark}
  
  \begin{theorem} \label{teo_criterio_de_regularidade}
    Let $(u,v)$ be a weak solution of \eqref{problema_P_controlado} with $f \in L^q(Q)$, $q > 5/2$. If, additionally, $u^0,v^0 \in W^{2 - 2/q,q}$ and $u^s \in L^q(Q)$ then $v \in X_q$ and $u \in L^{\infty}(Q)$. This implies, in particular, that $\nabla v \in L^{5q/(5-q)}(Q) \hookrightarrow L^5(Q)$, $u \in X_q$ and that $(u,v)$ is the unique strong solution of \eqref{problema_P_controlado}. Moreover, there exists $C = C(\norma{u^s}{L^q(Q)},\norma{f}{L^q(Q)},\norma{\nabla u}{L^{5/4}(Q)}) > 0$, which is continuous and increasing with respect to each entry, $\norma{u^s}{L^q(Q)}$, $\norma{f}{L^q}$ and $\norma{\nabla u}{L^{5/4}(Q)}$, such that
    \begin{equation} \label{estimativa_solucao_forte_em_relacao_a_f_nabla_u}
      \norma{(u,v)}{X_q \times X_q} \leq C(\norma{u^s}{L^q(Q)},\norma{f}{L^q(Q)},\norma{\nabla u}{L^{5/4}(Q)}).
    \end{equation}
  \end{theorem}
  \begin{proof}[\bf Proof]
    Analogously to Lemma \ref{lema_maquina_bootstrap}, we are going to prove that $(u,v) \in X_q \times X_q$ and, since the number of steps of the procedure of gaining regularity will be finite, the proof of \eqref{estimativa_solucao_forte_em_relacao_a_f_nabla_u} is a consequence of the estimates furnished by Lemmas \ref{lema_regularidade_eq_calor}, \ref{lema_regularidade_w_in_X_p} and \ref{lema_maquina_bootstrap}.
    
    Considering the regularity $v \in L^{\infty}(Q)$ given by the regularity of the weak solution $(u,v)$ of \eqref{problema_P_controlado}, since by hypothesis $u^s, f \in L^q(Q)$ then, by applying Lemma \ref{lema_regularidade_eq_calor} to the $v$-equation of \eqref{problema_P_controlado}, we conclude that
    \begin{equation} \label{regularidade_v_X_q}
      v \in X_q.
    \end{equation}
    Now denote $p_0 = sq > 5/2$. With $u \in L^{p_0}(Q)$ and \eqref{regularidade_v_X_q} we can apply Lemma \ref{lema_maquina_bootstrap} to conclude that the $u$-equation of \eqref{problema_P_controlado} is satisfied $a.e.$ $(t,x) \in Q$ and $u \in X_{q p_0/(q + p_0)}$.
    
    At least in this first iterations, we assume that we are in the case in which we have $q p_0/(q + p_0) < 5/2$. Then, if now we apply Lemma \ref{lema_regularidade_w_in_X_p} we obtain $u \in L^r(Q),$ with $r = \frac{5q}{5q + 5 p_0 - 2 q p_0} \ p_0.$ Since $q > 5/2$ we can say that that there is $\alpha > 1$ such that $q = 5\alpha/2$. Then using the fact that $p_0 \geq q$, we obtain
    \begin{align*}
      r & = \frac{5q}{5q - 5(\alpha - 1)p_0} \ p_0 = \frac{q}{q - (\alpha - 1)p_0} \ p_0 = \left ( 1 + \frac{(\alpha - 1)p_0}{q - (\alpha - 1)p_0} \right ) \ p_0 \\
      & > \left ( 1 + \frac{(\alpha - 1)p_0}{q} \right ) \ p_0 \geq (1 + \alpha - 1) p_0  = \alpha p_0
    \end{align*}
    Define $p_1 = \alpha p_0$. Since $r \geq p_1$ we have, in particular $u \in L^{p_1}(Q)$. Proceeding by induction, if we have $u \in L^{p_{n-1}}(Q)$, with $p_{n-1} = \alpha^{n-1} p_0 \geq q$ satisfying ${q p_{n-1}/(q + p_{n-1})} \leq 5/2$ then we can apply Lemmas \ref{lema_maquina_bootstrap} and \ref{lema_regularidade_w_in_X_p} and conclude that $u \in L^{p_n}(Q), \mbox{ with } p_n = \alpha^n p_0$.
    
    As a consequence, we can apply Lemma \ref{lema_maquina_bootstrap} and obtain $u \in X_{q p_n/(q + p_n)}$. Since $\alpha > 1$, we conclude that $p_n=\alpha^n p_0 = \alpha^n s q$ grows as $n$ increases in such a way that there is an index $n_0$ such that
      ${q p_{n_0-1}/(q + p_{n_0-1})} \leq 5/2$
    but applying the result proved by induction we conclude that
    \begin{equation*}
      u \in X_{q p_{n_0}/(q + p_{n_0})}, \mbox{ with } \frac{q p_{n_0}}{(q + p_{n_0})} > 5/2.
    \end{equation*}
    Hence, applying Lemma \ref{lema_regularidade_w_in_X_p} we obtain $u \in L^{\infty}(Q)$.
    
    Finally, once we have $u \in L^{\infty}(Q)$, we use \eqref{regularidade_v_X_q} and Lemma \ref{lema_maquina_bootstrap} to conclude that $u \in X_q$, finishing the proof.
  \end{proof}

  Now, to prove Theorem \ref{teo_dependencia_em_relacao_ao_controle}, it suffices to eliminate the dependence on $\norma{\nabla u}{L^{5/4}(Q)}$ in \eqref{estimativa_solucao_forte_em_relacao_a_f_nabla_u}. To this end, we introduce the auxiliary variable $z = \sqrt{v + \alpha^2}$, for some $\alpha > 0$, and we consider the (regularized) ``free-energy'' functional 
  \begin{equation*}
    E(u,z)(t) = \frac{s}{4} \D{\int_{\Omega}}{g(u(t,x)) \ dx} + \frac{1}{2} \int_{\Omega}{\norm{\nabla z(t,x)}{}^2} \ dx,
  \end{equation*}
  where
  \begin{equation*}
    g(u) = \left \{ 
    \begin{array}{ll}
      (u+1)ln(u+1) - u, & \mbox{if } s = 1,  \\
      \dfrac{u^s}{s(s-1)}, & \mbox{if } s > 1.
    \end{array}
    \right.
  \end{equation*}
  
  \begin{lemma} \label{lema_desigualdade_energia_solucao_forte}
    Assume \eqref{assumptions_set} and the existence of the strong solution $(u,v)$ of \eqref{problema_P_controlado} associated to $f \in L^q(Q)$. Let $z = \sqrt{v + \alpha^2}$, for some $\alpha > 0$ small enough. Then we have that $z$ satisfies
    \begin{equation} \label{limitacao_z_L_infty}
      0 < \alpha \leq z(t,x) \leq \mathcal{K}_1(\norma{f}{L^q(Q)},\norma{v^0}{W^{2-2/q,q}})
    \end{equation}
    and, moreover, there is $\alpha_0 > 0$, independent of $(u,v,f)$, such that if $0 < \alpha \leq \alpha_0$ then $(u,z)$ satisfies the energy inequality
    \begin{equation} \label{desigualdade_energia_solucao_forte}
      \begin{array}{r}
        E(u,z)(t_2) + \beta \D{\int_{t_1}^{t_2} \int_{\Omega}}{\norm{\nabla [u + 1]^{s/2}}{}^2 \ dx} \ dt + \dfrac{1}{4} \D{\int_{t_1}^{t_2} \int_{\Omega}}{u^s \norm{\nabla z}{}^2 \ dx} \ dt \\[6pt]
        + \beta \Big ( \D{\int_{t_1}^{t_2} \int_{\Omega}}{\norm{D^2 z}{}^2 \ dx} \ dt + \D{\int_{t_1}^{t_2} \int_{\Omega}}{\frac{\norm{\nabla z}{}^4}{z^2} \ dx} \ dt \Big ) \\[6pt]
        \leq E(u,z)(t_1) + \mathcal{K}_1(\norma{f}{L^q(Q)},\norma{v^0}{W^{2-2/q,q}}),
      \end{array}
    \end{equation}
    for $a.e.$ $t_1,t_2 \in [0,T]$, with $t_2 > t_1$; where $\mathcal{K}_1(\norma{f}{L^q(Q)},\norma{v^0}{W^{2-2/q,q}})$ is a continuous and increasing function with respect to $\norma{f}{L^q(Q)}$ and $\beta > 0$ is a constant, independent of $(u,v,f)$.
  \end{lemma}
  \begin{proof}[\bf Sketch of the proof]
    From \cite{guillen2022optimal}, given $f \in L^q(Q)$, one has the existence of weak solutions of \eqref{problema_P_controlado} satisfying the pointwise bound \eqref{limitacao_z_L_infty} and the energy inequality \eqref{desigualdade_energia_solucao_forte}. In the present lemma, we state that the unique strong solution  $(u,v)$ of \eqref{problema_P_controlado} given in Theorem~\ref{teo_criterio_de_regularidade} also satisfies \eqref{limitacao_z_L_infty} and \eqref{desigualdade_energia_solucao_forte}. Indeed, due to the strong regularity of $(u,v)$, the ideas of \cite{guillen2022optimal} can be applied directly, yielding the desired result. Next we give a sketch of the proof.

    First, using the variable $z = \sqrt{v^2 + \alpha^2}$, we rewrite  \eqref{problema_P_controlado} 
    as the following equivalent $(u,z)$ controlled problem
    \begin{equation} \label{problema_P_u_z_controlado}
      \left \{ \begin{array}{l}
        \partial_t u - \Delta u = - \nabla \cdot (u \nabla (z)^2) \\
        \partial_t z - \Delta z - \dfrac{\norm{\nabla z}{}^2}{z}
         = - \dfrac{1}{2} u^s \left( z - \dfrac{\alpha^2}{z} \right) + \dfrac{1}{2} f \left(z - \dfrac{\alpha^2}{z} \right) 1_{\Omega_c} \\
        \partial_\eta u  |_{\Gamma} = \partial_\eta z  |_{\Gamma} = 0, \ u(0) = u^0, z(0) = \sqrt{v^0 + \alpha^2}.
      \end{array} \right.
    \end{equation}
    
    To prove \eqref{limitacao_z_L_infty} one can adapt the result of \cite{guillen2022optimal}, where the pointwise estimates \eqref{limitacao_z_L_infty} are proved for the function $z_m$ which appears there, to the function $z$ of \eqref{problema_P_u_z_controlado}.
 
    For brevity, we just give an idea of the proof of \eqref{desigualdade_energia_solucao_forte} for the case $s = 1$. We begin by testing the $z$-equation of \eqref{problema_P_u_z_controlado} by $- \Delta z$, and using integration by parts, Holder's inequality and \eqref{limitacao_z_L_infty} in the right hand side terms, one has 
      \begin{align*}
       & \frac{1}{2} \frac{d}{dt} \norma{\nabla z}{L^2}^2 + \dfrac{1}{2} \norma{\Delta z}{L^2}^2 + \int_{\Omega}{\frac{\norm{\nabla z}{}^2}{z} \Delta z \ dx} + \frac{1}{2} \int_{\Omega}{u \norm{\nabla z}{}^2 \ dx} \\
       & \leq \frac{1}{4} \int_{\Omega}{\nabla (z^2) \cdot \nabla u \ dx} + \frac{1}{2} \alpha \int_{\Omega}{\norm{\nabla z}{} \norm{\nabla u}{} \ dx} + \mathcal{K}_1^2 \norma{f}{L^2}^2.
      \end{align*}
      Then, applying Lemma \ref{lema_termo_fonte_final}, we obtain
      \begin{equation} \label{estimativa_z_1_s=1}
        \begin{array}{l}
          \dfrac{1}{2} \dfrac{d}{dt} \norma{\nabla z}{L^2}^2 + C_1 \Big ( \D{\int_{\Omega}}{\norm{D^2 z}{}^2 \ dx} + \D{\int_{\Omega}}{\frac{\norm{\nabla z}{}^4}{z^2} \ dx} \Big ) + \frac{1}{2} \int_{\Omega}{u \norm{\nabla z}{}^2 \ dx} \\[6pt]
          \leq \dfrac{1}{4} \D{\int}_{\Omega}{\nabla (z^2) \cdot \nabla u \ dx} + \dfrac{1}{2} \alpha \D{\int}_{\Omega}{\norm{\nabla z}{} \norm{\nabla u}{} \ dx} + \mathcal{K}_1^2 \norma{f}{L^2}^2 + C_2 \norma{\nabla z}{L^2}^2.
        \end{array}
      \end{equation}

      Now we test the $u$-equation of \eqref{problema_P_u_z_controlado} by $ln(u + 1)$ and use 
      estimate \eqref{limitacao_z_L_infty} and 
      that $\dfrac{1}{u + 1} \leq \dfrac{1}{\sqrt{u + 1}}$, 
        \begin{align*}
          \frac{d}{dt} \int_{\Omega}{g(u) \ dx} & + C \int_{\Omega}{\norm{\nabla [u + 1]^{1/2}}{}^2 \ dx} = 2 \prodl{\frac{u + 1 - 1}{u + 1} z \nabla z}{\nabla u} \\
          & = \prodl{\nabla (z^2)}{\nabla u} - 2 \prodl{z \nabla z}{\frac{\nabla u}{u + 1}} \\
          & \leq \prodl{\nabla (z^2)}{\nabla u} + 2 \,
          \mathcal{K}_1(\norma{f}{L^q(Q)},\norma{v^0}{W^{2-2/q,q}}) 
          \norma{\nabla z}{L^2} \norma{\nabla [u + 1]^{1/2}}{L^2},
        \end{align*}
        Then, using adequately Young's inequality, we obtain
        \begin{equation}
          \frac{d}{dt} \int_{\Omega}{g(u) \ dx} + C \int_{\Omega}{\norm{\nabla [u + 1]^{1/2}}{}^2 \ dx} \leq \prodl{\nabla (z^2)}{\nabla u} + C \norma{\nabla z}{L^2}^2
          \label{estimativa_u_1_s=1}
        \end{equation}
        If we add inequality \eqref{estimativa_z_1_s=1} to $1/4$ times \eqref{estimativa_u_1_s=1}, then the terms $\D{\int_{\Omega}{\nabla (z^2) \cdot \nabla u \ dx}}$ cancel and we obtain
        \begin{align}
          & \nonumber \frac{d}{dt} E(u,z) + C \int_{\Omega}{\norm{\nabla [u + 1]^{1/2}}{}^2 \ dx} + \frac{1}{2} \int_{\Omega}{u \norm{\nabla z}{}^2 \ dx} \\
          & + C_1 \Big ( \int_{\Omega}{\norm{D^2 z}{}^2 \ dx} + \int_{\Omega}{\frac{\norm{\nabla z}{}^4}{z^2} \ dx} \Big ) \leq \frac{\alpha}{2} \int_{\Omega}{\norm{\nabla z}{} \norm{\nabla u}{} \ dx} + C \norma{\nabla z}{L^2}^2 + K_1^2 \norma{f}{L^2}^2 \label{aux_inequality_u_z_s=1} \\
          & \nonumber \leq \int_{\Omega}{\alpha \norm{\nabla [u + 1]^{1/2}}{} \norm{\sqrt{u + 1}}{} \norm{\nabla z}{} \ dx} + C \norma{\nabla z}{L^2}^2 + K_1^2 \norma{f}{L^2}^2.
        \end{align}
        We can deal with the first term in the right hand side of the inequality using Hölder's and Young's inequality,
        \begin{align*}
       &   \int_{\Omega}{{\alpha} \norm{\nabla [u + 1]^{1/2}}{} \norm{\sqrt{u + 1}}{} \norm{\nabla z}{} \ dx}
          \leq 
           \, \delta \int_{\Omega} \norm{\nabla [u + 1]^{1/2}}{}^2 \ dx
          + \alpha^2 \, C(\delta) \Big ( \int_{\Omega}{u \norm{\nabla z}{}^2 \ dx} + \int_{\Omega}{\norm{\nabla z}{}^2 \ dx} \Big ).
        \end{align*}
        Therefore, we can first choose $\delta > 0$ and then $\alpha > 0$ sufficiently small in order to use the terms on the left hand side of inequality \eqref{aux_inequality_u_z_s=1} to absorb the first two terms on the right hand side of the above inequality. Integrating the resulting inequality with respect to the variable $t$ from $t_1$ to $t_2$ gives us \eqref{desigualdade_energia_solucao_forte}, for the case $s = 1$.

        For the general case $s \geq 1$ we refer the reader to \cite{ViannaGuillen2023uniform} for details on the derivation of an energy inequality for the chemotaxis-consumption models \eqref{problema_P} and to \cite{guillen2022optimal} for the adaptation of this energy inequality to the controlled model \eqref{problema_P_controlado}.
  \end{proof}

  Using the results developed in the present section we are finally in position of proving Theorem~\ref{teo_dependencia_em_relacao_ao_controle}.
  
  \subsection{Proof of Theorem \ref{teo_dependencia_em_relacao_ao_controle}}

    Let $(u,v)$ be a weak solution of \eqref{problema_P_controlado} with $f,u^s \in L^q(Q)$, $q > 5/2$. From 
    Theorem~\ref{teo_criterio_de_regularidade},  $(u,v) \in X_q \times X_q$ is the strong solution of \eqref{problema_P_controlado} satisfying the intermediate estimate \eqref{estimativa_solucao_forte_em_relacao_a_f_nabla_u}.  Now, to prove the final estimate  \eqref{estimativa_solucao_forte_em_relacao_a_f}, it suffices to prove that $\norma{\nabla u}{L^{5/4}(Q)}$ can be estimated in terms of $\norma{f}{L^q(Q)}$. We analyze the cases $s \in [1,2)$ and $s \geq 2$ separately.
    
    \hspace{12pt}
    
    \noindent {\bf Case $\boldsymbol{s \in [1,2)}$:} Integrating the $u$-equation of \eqref{problema_P_controlado}, and using that $u\ge 0$, we obtain
    \begin{equation} \label{estimativa_L_infty_L_2_u_1/2}
      \norma{u^{1/2}}{L^{\infty}(L^2)}^2 = \int_{\Omega}{u^0(x) \ dx}.
    \end{equation}
    From \eqref{desigualdade_energia_solucao_forte} with $t_1 = 0$ and $t_2 = t$, we have, in particular,
    \begin{equation} \label{gradiente_u^s/2_L^2}
      \begin{array}{c}
        \norma{(u + 1)^{s/2}}{L^{\infty}(L^2)}^2 \leq E(u,z)(0) + \mathcal{K}_1(\norma{f}{L^q(Q)}, \norma{v^0}{W^{2-2/q,q}}), \quad\mbox{for } s > 1, \\[6pt]
        \norma{\nabla [u + 1]^{s/2}}{L^2(Q)}^2 \leq E(u,z)(0) + \mathcal{K}_1(\norma{f}{L^q(Q)}, \norma{v^0}{W^{2-2/q,q}}), \quad\mbox{for } s \geq 1.
     \end{array}
    \end{equation}
    Estimate \eqref{gradiente_u^s/2_L^2} (jointly with \eqref{estimativa_L_infty_L_2_u_1/2} for the case $s = 1$) implies that there exists $C= C(\norma{f}{L^q(Q)}) > 0$ ($C$ also depends on $(u^0,v^0)$, but since the initial data are fixed we omit it from now on) which is continuous and increasing with respect to $\norma{f}{L^q(Q)}$ and such that
    \begin{equation*}
      \norma{(u + 1)^{s/2}}{L^{\infty}(L^2) \cap L^2(H^1)} \leq C(\norma{f}{L^q(Q)}).
    \end{equation*}
    Using the latter and the interpolation inequality \eqref{interpolacao_norma_10/3}, we have
    \begin{equation} \label{u_L^5s/3}
      \norma{(u + 1)}{L^{5s/3}(Q)} \leq C(\norma{f}{L^q(Q)}).
    \end{equation}
    Now, recalling that $s \in [1,2)$, we consider the relation
    \begin{equation} \label{identidade_nabla_u}
      \nabla u = \nabla (u + 1) = \nabla \big ( (u + 1)^{s/2} \big )^{2/s} = \dfrac{2}{s} (u + 1)^{1 - s/2} \ \nabla (u + 1)^{s/2}.
    \end{equation}
    From \eqref{u_L^5s/3} we have that there is a constant ${C} ={C} (\norma{f}{L^q(Q)}) > 0$ which is continuous and increasing with respect to $\norma{f}{L^q(Q)}$ such that
    \begin{equation} \label{u^1-s/2}
      \norma{(u + 1)^{1-s/2}}{L^{10s/(6-3s)}(Q)} \leq {C}(\norma{f}{L^q(Q)}).
    \end{equation}
    From \eqref{identidade_nabla_u}, $\nabla u$ can be written as the product of $(u + 1)^{1-s/2}$ and $\nabla u$. By Holder's inequality, 
    \begin{equation*}
      \norma{\nabla u}{L^{5s/(3+s)}(Q)} \leq \norma{(u + 1)^{1-s/2}}{L^{10s/(6-3s)}(Q)} \norma{\nabla (u + 1)^{s/2}}{L^2(Q)}.
    \end{equation*}
    Hence, using \eqref{u^1-s/2} and \eqref{gradiente_u^s/2_L^2} in the above inequality, we conclude that there exists a constant $C = C(\norma{f}{L^q(Q)}) > 0$ which is continuous and increasing with respect to $\norma{f}{L^q(Q)}$ satisfying
    \begin{equation*}
      \norma{\nabla u}{L^{5s/(3 + s)}(Q)} \leq C(\norma{f}{L^q(Q)}).
    \end{equation*}
    Since $s \geq 1$, we have $5s/(3+s) \geq 5/4$ and this implies,
    \begin{equation} \label{dependencia_nabla_u_em_f_s<2}
      \norma{\nabla u}{L^{5/4}(Q)} \leq C(\norma{f}{L^q(Q)}).
    \end{equation}
    Therefore, using \eqref{dependencia_nabla_u_em_f_s<2} in \eqref{estimativa_solucao_forte_em_relacao_a_f_nabla_u} we conclude \eqref{estimativa_solucao_forte_em_relacao_a_f}.
    
    \hspace{12pt}
    
    \noindent {\bf Case $\boldsymbol{s \geq 2}$:} From \eqref{limitacao_z_L_infty} and \eqref{desigualdade_energia_solucao_forte} with $t_1 = 0$, there exists $C_1 = C_1(\norma{f}{L^q(Q)}) > 0$ such that
    \begin{equation} \label{estimativa_grad_z_e_termo_misto}
      \norma{z}{L^{\infty}(Q)}, \norma{\nabla z}{L^{\infty}(L^2)}, \norma{u^{s/2} \nabla z}{L^2(Q)} \leq C_1(\norma{f}{L^q(Q)}).
    \end{equation}
    Now, let us consider the sets
    \begin{equation*}
      \{ 0 \leq u \leq 1 \} = \Big \{ (t,x) \in Q \ \Big | \ 0 \leq u(t,x) \leq 1 \Big \}
    \end{equation*}
    and $\{ u \geq 1 \} = \Big \{ (t,x) \in Q \ \Big | \ u(t,x) \geq 1 \Big \}$. Note that, since $s \geq 2$, we have
    \begin{align*}
      & \int_0^T \int_{\Omega} u(t,x)^2 \norm{\nabla z(t,x)}{}^2 dx \ dt \\
      & \leq \int_0^T \int_{\{ 0 \leq u \leq 1 \}}{\norm{\nabla z(t,x)}{}^2 dx \ dt} + \int_0^T \int_{\{ u \geq 1 \}}{u(t,x)^s \norm{\nabla z(t,x)}{}^2 dx \ dt} \\
      & \leq \int_0^T \int_{\Omega}{\norm{\nabla z(t,x)}{}^2 \ dx} + \int_{\Omega}{u(t,x)^s \norm{\nabla z(t,x)}{}^2 \ dx}.
    \end{align*}
    Thus, by \eqref{estimativa_grad_z_e_termo_misto} we conclude that there is $C > 0$ such that
    \begin{equation} \label{limitacao_termo_misto_u_v_s_geq_2}
      \norma{u \nabla z}{L^2(Q)} \leq C \ C_1(\norma{f}{L^q(Q)}).
    \end{equation}
    
    Now we test the $u$-equation of \eqref{problema_P_controlado} by $ u$ and obtain
    \begin{align*}
      \dfrac{1}{2} \dfrac{d}{dt} \norma{u}{L^2 }^2 + \norma{\nabla u}{L^2 }^2 = 2 \int_{\Omega}{u z \nabla z \cdot \nabla u \ dx} \leq C \norma{z}{L^{\infty}(Q)}^2 \int_{\Omega}{u^2 \norm{\nabla z}{}^2 dx} + \frac{1}{2} \norma{\nabla u}{L^2 }^2.
    \end{align*}
    Hence we have
    \begin{equation*}
      \dfrac{d}{dt} \norma{u}{L^2 }^2 + \norma{\nabla u}{L^2 }^2 \leq C \norma{z}{L^{\infty}(Q)}^2 \int_{\Omega}{u^2 \norm{\nabla v}{}^2 dx}.
    \end{equation*}
    Integrating with respect to $t$, we conclude from \eqref{estimativa_grad_z_e_termo_misto} that there is $C_2 = C_2(\norma{f}{L^q(Q)}) > 0$ such that
      $\norma{\nabla u}{L^2(Q)} \leq C_2(\norma{f}{L^q(Q)}).$
    This implies, in particular, that we have \eqref{dependencia_nabla_u_em_f_s<2} also for the case $s \geq 2$ and therefore, using again \eqref{dependencia_nabla_u_em_f_s<2} in \eqref{estimativa_solucao_forte_em_relacao_a_f_nabla_u} leads us to \eqref{estimativa_solucao_forte_em_relacao_a_f}.

 %
  \section{Proof of Theorem \ref{teo_existencia_controle_otimo}}
  \label{section: existencia controle otimo}
    From \eqref{hipotese_conjunto_admissivel_nao_vazio} and since the functional $J$ in \eqref{problema_de_minimizacao} is nonnegative, then
      $J_{inf} : = \D{\inf_{(u,v,f) \in S_{ad}}} J(u,v,f) \geq 0$
    is well defined and there is a minimizing sequence $\{ (u_n,v_n,f_n) \} \subset S_{ad}$ satisfying
    \begin{equation} \label{problema_P_controlado_sequencia_minimizante}
      \left\{\begin{array}{l}
        \partial_t u_n - \Delta u_n  = - \nabla \cdot (u_n \nabla v_n), \quad
        \partial_t v_n - \Delta v_n  = - u_n^s v_n + f_n v_n 1_{\Omega_c}, \\
        \partial_{\boldsymbol{n}} u_n |_{\Gamma}  =  \partial_{\boldsymbol{n}} v_n |_{\Gamma} = 0, \quad
        u_n(0)  = u^0, \quad v_n(0) = v^0,
      \end{array}\right.
    \end{equation}
    and $\lim_{n \to \infty}{J(u_n,v_n,f_n)} = J_{inf}.$ Next we prove that there is $(\overline{u},\overline{v},\overline{f}) \in S_{ad}$, that will be defined as the limit of a subsequence of $\{ (u_n,v_n,f_n) \}_n$, such that $J(\overline{u},\overline{v},\overline{f}) = J_{inf}$.
    
    In fact, from the definition of $J$ and the hypothesis \eqref{hipotese_sobre_os_custos_gamma}, we have
    \begin{equation*}
      \{ u_n^s \}_n \mbox{ is bounded in } L^q(Q),
    \end{equation*}
    \begin{equation} \label{limitacao_f_n}
      \{ f_n \}_n \mbox{ is bounded in } L^q(Q).
    \end{equation}
    Since $(u_n,v_n,f_n) \in S_{ad}$, $(u_n,v_n)$ is the strong solution of \eqref{problema_P_controlado} with control $f_n$. Then, from   \eqref{limitacao_f_n} and Theorem \ref{teo_dependencia_em_relacao_ao_controle} we obtain
    \begin{equation} \label{limitacao_u_n_v_n}
      \{ u_n \}_n \mbox{ and } \{ v_n \}_n \mbox{ are bounded in } X_q.
    \end{equation}
    We recall that since $\mathcal{F}$ is a closed and convex subset of $L^q(Q)$ then $\mathcal{F}$ is also weakly closed in $L^q(Q)$. Therefore, accounting for the $n$-independent bounds \eqref{limitacao_f_n} and \eqref{limitacao_u_n_v_n} we conclude that there exists $(\overline{u},\overline{v},\overline{f}) \in X_q \times X_q \times \mathcal{F}$ such that, up to a subsequence, we have the weak convergences as $n\to +\infty$:
    \begin{equation} \label{convergencia_sequencia_minimizante_fraca}
      \begin{array}{c}
        (u_n,v_n) \rightarrow (\overline{u},\overline{v}) \mbox{ weakly* in } L^{\infty}(W^{2 - 2/q,q}) \times L^{\infty}(W^{2 - 2/q,q}) \\
        (u_n,v_n) \rightarrow (\overline{u},\overline{v}) \mbox{ weakly in } L^q(W^{2,q}) \times L^q(W^{2,q}), \\
        (\partial_t u_n, \partial_t v_n) \rightarrow (\partial_t \overline{u},\partial_t \overline{v}) \mbox{ weakly in } L^q(Q) \times L^q(Q), \\
        f_n \rightarrow \overline{f} \mbox{ weakly in } L^q(Q).
      \end{array}
    \end{equation}
    Since $q > 5/2$, we have $W^{2 - 2/q,q}(\Omega)$ compactly embedded in $C^0(\overline{\Omega})$,  hence from the compactness Lemma \ref{lema_Simon}, we get 
    \begin{equation} \label{convergencia_sequencia_minimizante_forte_nas_continuas}
      (u_n,v_n) \rightarrow (\overline{u},\overline{v}) \mbox{ strongly in } C(\overline{Q}) \times C(\overline{Q}).
    \end{equation}
    Following the proofs of Lemmas \ref{lema_regularidade_w_in_X_p} and \ref{lema_regularidade_nabla_w_in_X_p}, one can conclude from \eqref{limitacao_u_n_v_n} that $\{ (u_n,v_n) \}_n$ is compactly embedded in $L^{5q/(5-2q)}(Q)$ and $\{ (\nabla u_n, \nabla v_n) \}_n$ is compactly embedded in $L^p(Q)$, for all $p < 5q/(5-q)$. Since $5q/(5-2q) > 5q/(5-q)$, we have
    \begin{equation} \label{injecao_compacta_u_n_v_n}
      \{ (u_n,v_n) \}_n \mbox{ is compactly embedded in } L^p(W^{1,p}) \times L^p(W^{1,p}), \ \forall p < 5q/(5-q).
    \end{equation}
    Then, from \eqref{limitacao_u_n_v_n}, \eqref{injecao_compacta_u_n_v_n} and Lemma \ref{lema_Simon} we obtain
    \begin{equation} \label{convergencia_sequencia_minimizante_forte_derivadas}
      (u_n, v_n) \rightarrow (\overline{u}, \overline{v}) \mbox{ strongly in } L^p(W^{1,p}) \times L^p(W^{1,p}), \ \forall p < 5q/(5-q).
    \end{equation}
    Recalling that $q > 5/2$, we have $5q/(5-q) > 2q$ and thus, from \eqref{convergencia_sequencia_minimizante_forte_derivadas} we also have, in particular,
    \begin{equation} \label{convergencia_sequencia_minimizante_forte_dos_gradientes}
      (\nabla u_n, \nabla v_n) \rightarrow (\nabla \overline{u}, \nabla \overline{v}) \mbox{ strongly in } L^{2q}(Q) \times L^{2q}(Q).
    \end{equation}
    From the above strong convergences we conclude that $\overline{u}(0) = u^0$ and $\overline{v}(0) = v^0$. Moreover, since
    \begin{equation*}
      0 = \int_{\Gamma} \partial_{\boldsymbol{n}} u_n(t, \cdot) \ \varphi |_{\Gamma} \ d\Gamma = \int_{\Omega} \Delta u_n(t,x) \varphi(x) \ dx + \int_{\Omega} \nabla u_n(t,x) \cdot \nabla \varphi(x) \ dx, \ \forall \varphi \in C^{\infty}_c(\Omega),
    \end{equation*}
    $a.e.$ $t \in (0,T)$, the weak convergence \eqref{convergencia_sequencia_minimizante_fraca} implies that $\partial_{\boldsymbol{n}} \overline{u} |_{\Gamma} = 0$. Analogously, we also have $\partial_{\boldsymbol{n}} \overline{v} |_{\Gamma} = 0$.
    
    With the convergences \eqref{convergencia_sequencia_minimizante_fraca}, \eqref{convergencia_sequencia_minimizante_forte_nas_continuas} and \eqref{convergencia_sequencia_minimizante_forte_dos_gradientes} we pass to the limit in the nonlinear terms of \eqref{problema_P_controlado_sequencia_minimizante} and prove that
    \begin{equation*}
      \begin{array}{c}
        \nabla u_n \cdot \nabla v_n + u_n \Delta v_n \rightarrow \nabla \overline{u} \cdot \nabla \overline{v} + \overline{u} \Delta \overline{v} \mbox{ weakly in } L^q(Q), \\
        u_n^s v_n \rightarrow \overline{u}^s \, \overline{v} \mbox{ strongly in } C(\overline{Q}), \\
        f_n v_n 1_{\Omega_c} \rightarrow \overline{f} \overline{v} 1_{\Omega_c} \mbox{ weakly in } L^q(Q).
      \end{array}
    \end{equation*}
  Since passing to the limit in the linear terms of \eqref{problema_P_controlado_sequencia_minimizante} is rather standard, we have proved that $(\overline{u},\overline{v}) \in X_q \times X_q$ is the strong solution of \eqref{problema_P_controlado} with control $\overline{f} \in \mathcal{F}$, that is, $(\overline{u},\overline{v}, \overline{f}) \in S_{ad}$. Hence, we have, in particular,
    \begin{equation} \label{desigualdade_limite_sequencia_minimizante_1}
      \inf_{(u,v,f) \in S_{ad}} J(u,v,f) \leq J(\overline{u},\overline{v},\overline{f}).
    \end{equation}
    On the other hand, using the fact that the functional $J$ is lower weakly semicontinuous, we also have
   $
      J(\overline{u},\overline{v},\overline{f}) \leq \D{\inf_{(u,v,f) \in S_{ad}}} J(u,v,f),
  $
    thus, jointly to     
     \eqref{desigualdade_limite_sequencia_minimizante_1}, one has that $(\overline{u},\overline{v},\overline{f})$ is a global optimum.


\section{First order necessary conditions for any local optimum}
  \label{section: multiplicadores de Lagrange}
  
   To derive the first order necessary optimality conditions for a local optimal solution $(\overline{u}, \overline{v}, \overline{f})$ of \eqref{problema_de_minimizacao}, we use a generic Lagrange Multipliers theorem given by \cite{zowe1979regularity}  that we introduce in Subsection \ref{subsec:abstract setting and Lagrange multipliers}. Then, in Subsection \ref{subsec: regular point} we prove that any local optimal solution is a regular point (see Definition \ref{defi_regular_point} below), in Subsection \ref{subsec: existencia_multiplicador_de_Lagrange_aplicado} we prove Theorem \ref{teo_existencia_multiplicador_de_Lagrange_aplicado}, and in Subsection \ref{subsec: regularidade_adicional_multiplicadores_de_Lagrange} we finally prove Theorem \ref{teo_regularidade_adicional_multiplicadores_de_Lagrange}.


\subsection{
 Lagrange multipliers theorem}
  \label{subsec:abstract setting and Lagrange multipliers}
  
  Let us consider the following abstract optimization problem:
  \begin{equation} \label{problema_de_minimizacao_abstrato}
    \min_{r \in \mathbb{M}} J(r) \ \mbox{ subject to } \ G(r) = 0,
  \end{equation}
  where $J: \mathbb{X} \rightarrow \mathbb{R}$ is a functional, $G: \mathbb{X} \rightarrow \mathbb{Y}$ is an operator $\mathbb{X}$ and $\mathbb{Y}$ are Banach spaces and $\mathbb{M} \subset \mathbb{X}$ is a closed and convex subset. Note that the admissible set for problem \eqref{problema_de_minimizacao_abstrato} is $S = \{ r \in \mathbb{M} \ | \ G(r) = 0 \}$.
  
  Next we define the Lagrangian functional, the Lagrange multipliers and the so called regular points.
  \begin{definition}{\bf (Lagrangian)}
    The functional $\mathcal{L}: \mathbb{X} \times \mathbb{Y}' \rightarrow \mathbb{R}$, given by
    \begin{equation} \label{Lagrangian_functional}
      \mathcal{L}(r,\xi) = J(r) - \dist{\xi}{G(r)}_{\mathbb{Y}'},
    \end{equation}
    is called the Lagrangian functional related to problem \eqref{problema_de_minimizacao_abstrato}.
    \hfill $\square$
  \end{definition}
  
  \begin{definition}{\bf (Lagrange multipliers)}
    Let $\overline{r} \in S$ be a local optimal solution of problem \eqref{problema_de_minimizacao_abstrato}. Suppose that $J$ and $G$ are Fréchet differentiable in $\overline{r}$, the derivatives being denoted by $J'(\overline{r})$ and $G'(\overline{r})$, respectively. Then, $\xi \in \mathbb{Y}'$ is called a Lagrange multiplier for \eqref{problema_de_minimizacao_abstrato} at the point $\overline{r}$ if
    \begin{equation} \label{expressao_multiplicador_de_Lagrange}
      \mathcal{L}'(\overline{r},\xi)[c] = J'(\overline{r})[c] - \dist{\xi}{G'(\overline{r})[c]}_{\mathbb{Y}'} \geq 0, \ \forall c \in \mathcal{C}(\overline{r}),
    \end{equation}
    where $\mathcal{C}(\overline{r}) = \{ \theta (r - \overline{r}) \ | \ r \in \mathbb{M}, \ \theta \geq 0 \}$ is the conical hull of $\overline{r} \in \mathbb{M}$.
    \hfill $\square$
  \end{definition}
  
  \begin{definition}{\bf (Regular point)} \label{defi_regular_point}
    A point $\overline{r} \in \mathbb{M}$ is called a regular point if $G'(\overline{r})[\mathcal{C}(\overline{r})] = \mathbb{Y}$.
    \hfill $\square$
  \end{definition}
  Finally, we state the theorem on the existence of Lagrange multipliers.
  \begin{theorem}{\bf (\cite{zowe1979regularity})} \label{teo_existencia_multiplicadores_de_Lagrange}
    Let $\overline{r} \in S$ be a local optimal solution of problem \eqref{problema_de_minimizacao_abstrato}. Suppose that $J$ is Fréchet differentiable and $G$ is continuously Fréchet differentiable. If $\overline{r}$ is a regular point, then there exists a  Lagrange multiplier for problem \eqref{problema_de_minimizacao_abstrato} at $\overline{r}$.
  \end{theorem}


\subsection{Local optimal solutions are regular points}
  \label{subsec: regular point}
  
  To apply the theory of Subsection \ref{subsec:abstract setting and Lagrange multipliers} to our optimal control problem \eqref{problema_de_minimizacao}   
  and derive the first order necessary conditions for a local optimal solution of \eqref{problema_de_minimizacao}, we will reformulate \eqref{problema_de_minimizacao} using the abstract setting of \eqref{problema_de_minimizacao_abstrato}. Since we want $\mathbb{X}$ and $\mathbb{Y}$ to be Banach spaces, let us define them as $\mathbb{X} = \widetilde{X}_q \times \widetilde{X}_q \times L^q(Q), \ \mathbb{Y} = L^q(Q) \times L^q(Q)$, where $\widetilde{X}_q = \{ w \in X_q \ | \ \partial_{\boldsymbol{n}} w |_{\Gamma} = 0 \}$. Next we define the operator $G = (G_1,G_2): \mathbb{X} \rightarrow \mathbb{Y}$, where
    $G_1: \mathbb{X} \rightarrow L^q(Q)$ and $G_2: \mathbb{X} \rightarrow L^q(Q)$
  are defined for each $r = (u,v,f) \in \mathbb{X}$ as
  \begin{equation*}
    \left \{
    \begin{array}{rl}
      G_1(r) & \hspace{-2mm} = \partial_t u - \Delta u + \nabla \cdot (u \nabla v) \\
      G_2(r) & \hspace{-2mm} = \partial_t v - \Delta v + u^sv - fv \ 1_{\Omega_c}.
    \end{array}
    \right.
  \end{equation*}
  By using hypothesis $S_{ad}\not=\emptyset$, there exists $(\hat{u},\hat{v},\hat{f})\in S_{ad}$. Then, 
  we introduce the space 
  $$\widehat{X}_q = \{ w \in \widetilde{X}_q \ | \ w(0,x) = 0 \}$$ and we define $\mathbb{M}$, the closed and convex subset of $\mathbb{X}$, as $\mathbb{M} = (\widehat{u}, \widehat{v}, \widehat{f}) + \widehat{X}_q \times \widehat{X}_q \times (\mathcal{F} - \widehat{f}).$
 Therefore,  we rewrite the optimal control problem \eqref{problema_de_minimizacao} as
  \begin{equation} \label{problema_de_minimizacao_abstrato_2}
    \min_{r \in \mathbb{M}} J(r) \ \mbox{ subject to } \ G(r) = 0.
  \end{equation}
  
  We have the 
  differenciability of the functional $J$ and the operator $G$.
  \begin{lemma} \label{lema_J_diferenciavel}
    The functional $J: \mathbb{X} \rightarrow \mathbb{R}$ is Fréchet differentiable and the Fréchet derivative of $J$ in $\overline{r} = (\overline{u},\overline{v},\overline{f}) \in \mathbb{X}$ in the direction $c = (U,V,F) \in \mathbb{X}$ is 
    \begin{equation} \label{derivada_Frechet_J}
        J'(\overline{r})[c] 
        = \D{\int_0^T \int_{\Omega}} (g_\lambda\, U + g_\eta\, V) dx \ dt
         + \gamma_f \D{\int_0^T \int_{\Omega_c}} sgn(\overline{f}) \norm{\overline{f}}{}^{q-1} F \ dx \ dt,
    \end{equation}
    where $g_\lambda,g_\eta$ are defined in \eqref{rhs-adjoint}.
  \end{lemma}
  \begin{proof}[\bf Proof]
    The functional $J$ is the sum of functionals of the form $A: L^p(Q) \rightarrow \mathbb{R}$ given by
    \begin{equation*}
      A(w) = \int_0^T \int_{\Omega} \norm{w(t,x)}{}^p \ dx \ dt.
    \end{equation*}
    Hence, the Fréchet differentiability of the functional $J$ from $L^{sq}(Q) \times L^2(Q) \times L^q(Q)$ into $\mathbb{R}$ and the expression \eqref{derivada_Frechet_J} follow from the fact that $A$ is Fréchet differentiable from $L^p(Q)$ into $\mathbb{R}$ for all $p \in (1,\infty)$ with
    \begin{equation*}
      A'(w) h = p \int_0^T \int_{\Omega} sgn(w(t,x)) \norm{w(t,x)}{}^{p-1} \ h(t,x) \ dx \ dt, \quad \forall\, h \in L^p(Q).
    \end{equation*}
    This can be proved by applying the theory of Nemytskii operators (or superposition operators) in $L^p$-spaces given in \cite{goldberg1992nemytskij} and summarized in \cite[Subsection 4.3.3]{troltzsch2010optimal}. The Fréchet differentiability of $J$ from $\mathbb{X}$ into $\mathbb{R}$ follows from the fact that $\mathbb{X} \hookrightarrow L^{sq}(Q) \times L^2(Q) \times L^q(Q)$.
  \end{proof}
  
  \begin{lemma} \label{lema_G_diferenciavel}
   The operator $G: \mathbb{X} \rightarrow \mathbb{Y}$ is continuously Fréchet differentiable and the Fréchet derivative of $G$ in ${r} = ({u},{v},{f}) \in \mathbb{X}$ in the direction $c = (U,V,F) \in \mathbb{X}$ is the operator $G'({r})[c] = (G_1'({r})[c],G_2'({r})[c])$ given by
   \begin{equation} \label{derivada_Frechet_G}
    \left \{
    \begin{array}{rl}
      G_1'({r})[c] & \hspace{-2mm} = \partial_t U - \Delta U + \nabla \cdot (U \nabla{v}) + \nabla \cdot ({u} \nabla V) \\
      G_2'({r})[c] & \hspace{-2mm} = \partial_t V - \Delta V + s\, {u}^{s-1} U{v} +{u}^s V -{f} V \ 1_{\Omega_c} - F{v} \ 1_{\Omega_c}.
    \end{array}
    \right.
  \end{equation}
  \end{lemma}
  \begin{proof}[\bf Proof]
    Let us consider only the operator $G_1$. The proof of the Fréchet differentiability of the operator $G_2$ is analogous to the proof for $G_1$.
    
    We must show that
    \begin{equation} \label{propriedade_residuo_operador_G_1}
      \lim_{\norma{(U,V,F)}{\mathbb{X}} \to 0}{\frac{\norma{r_1((u,v,f),(U,V,F))}{L^q(Q)}}{\norma{(U,V,F)}{\mathbb{X}}}} = 0.
    \end{equation}
    where
      $r_1((u,v,f),(U,V,F)) = G_1(u + U,v + V,f + F) - G_1(u,v,f) - G_1'(u,v,f)[U,V,F]$.
    Expanding the terms of $r_1((u,v,f),(U,V,F))$ using the expressions of $G_1$ and $G_1'$ we get
    \begin{equation*}
      r_1((u,v,f),(U,V,F)) = \nabla \cdot (U \nabla V) = U \Delta V + \nabla U \cdot \nabla V.
    \end{equation*}
    Since $U,V \in X_q$, then $\Delta V \in L^q(Q)$ with
    \begin{equation} \label{Delta_V_em_L^q}
      \norma{\Delta V}{L^q(Q)} \leq C \norma{V}{X_q}.
    \end{equation}
    Since $q > 5/2$, from Lemma \ref{lema_regularidade_w_in_X_p}, we obtain $U \in L^{\infty}(Q)$ with
    \begin{equation} \label{U_em_L^infty}
      \norma{U}{L^{\infty}(Q)} \leq C \norma{U}{X_q}.
    \end{equation}
    Moreover, we have $L^{2q}(Q) \hookrightarrow L^{5q/(5-q)}(Q)$. Indeed, $5q/(5-q) > 5q/(5-5/2) = 2q$. This continuous injection jointly with Lemma \ref{lema_regularidade_nabla_w_in_X_p} gives us
    \begin{equation} \label{nabla_U_V_em_L^2q}
      \norma{\nabla U}{L^{2q}(Q)^3} \leq C \norma{U}{X_q} \quad \mbox{ and }
       \quad \norma{\nabla V}{L^{2q}(Q)^3} \leq C \norma{V}{X_q}.
    \end{equation}
    Hence, from \eqref{Delta_V_em_L^q}, \eqref{U_em_L^infty} and \eqref{nabla_U_V_em_L^2q} we conclude that
    \begin{equation*}
      \norma{r_1((u,v,f),(U,V,F))}{L^q(Q)} \leq C \norma{U}{L^{\infty}(Q)} \norma{\Delta V}{L^q(Q)} + C \norma{\nabla U}{L^{2q}(Q)^3} \norma{\nabla V}{L^{2q}(Q)^3}.
    \end{equation*}
    In particular, $\norma{r_1((u,v,f),(U,V,F))}{L^q(Q)} \leq C \,\norma{(U,V,F)}{\mathbb{X}}^2.$ Therefore we have
    \begin{equation*}
      0 \leq \frac{\norma{r_1((u,v,f),(U,V,F))}{L^q(Q)}}{\norma{(U,V,F)}{\mathbb{X}}} \leq C \,\norma{(U,V,F)}{\mathbb{X}},
    \end{equation*}
    for all $(U,V,F) \in \mathbb{X}$ such that $(u + U, v + V, f + F) \in \mathbb{X}$, which implies \eqref{propriedade_residuo_operador_G_1}.
  \end{proof}

  \begin{remark}
    To have \eqref{nabla_U_V_em_L^2q} it suffices that $q \geq 5/2$, but for \eqref{U_em_L^infty} it is crucial that $q > 5/2$.
    \hfill $\square$
  \end{remark}
  
  Next we prove the existence of Lagrange multipliers for the problem \eqref{problema_de_minimizacao_abstrato_2} associated to a local optimal solution $\overline{r} = (\overline{u},\overline{v},\overline{f}) \in S_{ad}$. Accounting for Lemmas \ref{lema_J_diferenciavel} and \ref{lema_G_diferenciavel} and Theorem \ref{teo_existencia_multiplicadores_de_Lagrange}, now it suffices to prove that $\overline{r}$ is a regular point, that is, 
  for each $(g_U,g_V) \in \mathbb{Y}$, there is $c = (U,V,F) \in \widehat{X}_q \times \widehat{X}_q \times \mathcal{C}(\overline{f})$ such that
  \begin{equation*}
    \left \{
    \begin{array}{rl}
      \partial_t U - \Delta U + \nabla \cdot (U \nabla \overline{v}) + \nabla \cdot (\overline{u} \nabla V) & \hspace{-2mm} = g_U \\
      \partial_t V - \Delta V + s \overline{u}^{s-1} U \overline{v} + \overline{u}^s V - \overline{f} V \ 1_{\Omega_c} - F \overline{v} \ 1_{\Omega_c} & \hspace{-2mm} = g_V.
    \end{array}
    \right.
  \end{equation*}
  where $\mathcal{C}(\overline{f}) = \{ \theta (f - \overline{f}) \ | \ f \in \mathcal{F}, \ \theta \geq 0 \}$ is the conical hull of $\overline{f} \in \mathcal{F}$. Since $0 \in \mathcal{C}(\overline{f})$, we can take $F = 0$ and therefore,  
  it suffices to prove that, for any $(g_U,g_V) \in \mathbb{Y}$, there is $ (U,V) \in \widehat{X}_q \times \widehat{X}_q $ such that
  \begin{equation} \label{problema_linearizado}
    \left \{
    \begin{array}{l}
      \partial_t U - \Delta U = - \nabla \cdot (U \nabla \overline{v}) - \nabla \cdot (\overline{u} \nabla V)  + g_U \\
      \partial_t V - \Delta V = - s \overline{u}^{s-1} U \overline{v} - \overline{u}^s V + \overline{f} V \ 1_{\Omega_c} + g_V.
    \end{array}
    \right.
  \end{equation}
  Problem \eqref{problema_linearizado} is called the linearized problem related to \eqref{problema_P_controlado}. Now we prove that $\overline{r}$ is a regular point. For this, we will use the general result Theorem \ref{teo_problema_linearizado_geral} given in 
  the Appendix \ref{appendixA}. Here, we consider the Banach space for weak solutions
  \begin{equation} \label{weak_solutions_space}
    W_2 = \{ v \in L^2(H^{1}) \ : \ \partial_t v \in L^{2}((H^1)') \}.
  \end{equation}

  \begin{remark}
    The space $W_2$ is continuously embedded in $C([0,T];L^2)$.
    \hfill $\square$
  \end{remark}
  
  \begin{theorem} \label{teo_ponto_regular}
    Let $\overline{r} = (\overline{u},\overline{v},\overline{f}) \in S_{ad}$.
    Then $\overline{r}$ is a regular point.
  \end{theorem}
  \begin{proof}[\bf Proof]
    Using Theorem \ref{teo_problema_linearizado_geral}, case $2a$, with $a_1 = b_1 = 0$, $\Vec{c}_1 = \nabla \overline{v} \in L^{5q/(5-q)}$, $d = \overline{u} \in L^{\infty}(Q)$, $a_2 = \overline{u}^s + \overline{f} 1_{\Omega_c} \in L^q(Q)$, $b_2 = s \overline{u}^{s-1} \, \overline{v} \in L^{\infty}(Q)$ and $\Vec{c}_2 = 0$,  we claim that there is a solution
    \begin{equation} \label{regularidade_solucao_forte_lema_problema_linearizado}
      (U,V) \in W_2 \times X_2
    \end{equation}
    of \eqref{problema_linearizado}. Therefore it suffices to prove that actually $(U,V) \in X_q \times X_q$. In fact, since $V \in X_2$, we have from Lemma \ref{lema_regularidade_w_in_X_p} that $V \in L^{10}(Q)$. Let $Z_1 = - s \overline{u}^{s-1} U \overline{v} - \overline{u}^s V + \overline{f} V \ 1_{\Omega_c} + g_V$ be the right hand side of the $V$-equation of \eqref{problema_linearizado}, then, accounting for the extra regularity of the coefficients (when compared to Theorem \ref{teo_problema_linearizado_geral}) we conclude that $Z_1 \in L^{10q/(10+q)}(Q)$ and, from Lemma \ref{lema_regularidade_eq_calor}, we have
    \begin{equation} \label{regularidade_V_forte_X_10q/10+q}
      V \in X_{10q/(10+q)}. 
    \end{equation}
    Note that $10q/(10 + q) < q$. We will enhance the regularity of $V$ and prove that $V \in X_q$ by induction. In fact suppose that $Z_1 \in L^{10q/(10 n + (5 - 4n)q)}(Q)$, with $10q/(10 n + (5 - 4n)q) < q$. From Lemma \ref{lema_regularidade_eq_calor} we have
    $
      V \in X_{10q/(10 n + (5 - 4n)q)}.
    $
    Using Lemma \ref{lema_regularidade_w_in_X_p} we have $V \in L^{10q/(10 n + (5 - 4(n+1))q)}(Q)$. Applying this regularity to the less regular term of $Z_1$, $\overline{f} V \ 1_{\Omega_c}$, we conclude that 
    $\overline{f} V \ 1_{\Omega_c} \in L^{10q/(10 (n+1) + (5 - 4(n+1))q)}(Q).$
     Thus, if $10q/(10 (n+1) + (5 - 4(n+1))q) < q$ then we conclude that $Z_1 \in L^{10q/(10 (n+1) + (5 - 4(n+1))q)}(Q)$.
    
    Therefore we have proved that, as long as $10q/(10n + (5 - 4n)q) < q$, if $Z_1 \in L^{10q/(10n + (5 - 4n)q)}(Q)$ then  $Z_1 \in L^{10q/(10(n+1) + (5 - 4(n+1))q)}(Q)$. Recalling that $q > 5/2$, if we study the function $n \mapsto 10q/(10n + (5 - 4n)q)$, we conclude that there exists $n_0,n_1 \in \mathbb{N}$ such that, $10q/(10n_0 + (5 - 4n_0)q) < q$ and $10q/(10n_1 + (5 - 4n_1)q) \geq q$. Thus we proved that the right hand side of the $V$-equation of \eqref{problema_linearizado} belongs to $L^q(Q)$. Finally, from Lemma \ref{lema_regularidade_eq_calor}, we have
    \begin{equation} \label{regularidade_V_X_q}
      V \in X_q.
    \end{equation}
    
    It remains to prove that $U \in X_q$. For this, we will analyze the right hand side of the $U$-equation of \eqref{problema_linearizado} and use \eqref{regularidade_solucao_forte_lema_problema_linearizado} and \eqref{regularidade_V_X_q}. The right hand side of the $U$-equation is
    \begin{equation*}
      g_U - U \Delta \overline{v} - \nabla U \cdot \nabla \overline{v} - \overline{u} \Delta V - \nabla \overline{u} \cdot \nabla V.
    \end{equation*}
    With the regularities obtained so far for $U$ and $V$, we have $g_U - \overline{u} \Delta V - \nabla \overline{u} \cdot \nabla V \in L^q(Q)$ and
    \begin{equation}
      Z_2 : = U \Delta \overline{v} + \nabla U \cdot \nabla \overline{v} \in L^{10q/(10 + q)}(Q).
    \end{equation}
    Again, we prove by induction that, as long as $10q/(10n + (5 - 4n)q) < p$, if $Z_2 \in L^{10q/(10n + (5 - 4n)q)}(Q)$ then we have $Z_2 \in L^{10q/(10(n+1) + (5 - 4(n+1))q)}(Q)$. Recalling that $q > 5/2$, if we study the function $n \mapsto 10q/(10n + (5 - 4n)q)$, we conclude that there exists $n_0,n_1 \in \mathbb{N}$ such that, $10q/(10n_0 + (5 - 4n_0)q) < q$ and $10q/(10n_1 + (5 - 4n_1)q) \geq q$ and thus we proved that the right hand side of the $U$-equation of \eqref{problema_linearizado} belongs to $L^q(Q)$. From Lemma \ref{lema_regularidade_eq_calor}, we conclude that $U \in X_q$. Finally, accounting for the linearity of \eqref{problema_linearizado}, one can deduce the uniqueness of the strong solution $(U,V) \in X_q \times X_q$.
  \end{proof}
  


  \subsection{Proof of Theorem \ref{teo_existencia_multiplicador_de_Lagrange_aplicado}}
\label{subsec: existencia_multiplicador_de_Lagrange_aplicado}

    The proof is divided in two steps: the existence of Lagrange multiplier and the uniqueness.

 \
    
    \noindent{\bf Step 1: Existence}

\vspace{3pt}
    
    From Lemmas \ref{lema_J_diferenciavel}, \ref{lema_G_diferenciavel} and Theorem \ref{teo_ponto_regular}, all the hypotheses of Theorem~\ref{teo_existencia_multiplicadores_de_Lagrange} are fulfilled. Therefore, there exists a Lagrange multiplier $\xi = (\lambda,\eta) \in L^{q'}(Q) \times L^{q'}(Q)$ satisfying, according to \eqref{expressao_multiplicador_de_Lagrange}, the inequality
    \begin{equation} \label{expressao_multiplicador_de_Lagrange_aplicada}
      \mathcal{L}'(\overline{r},\lambda,\eta)[c] = J'(\overline{r})[c] - \dist{\lambda}{G_1'(\overline{r})[c]}_{L^{q'}(Q)} - \dist{\eta}{G_2'(\overline{r})[c]}_{L^{q'}(Q)} \geq 0,
    \end{equation}
    for all $c = (U,V,F) \in \widehat{X}_q \times \widehat{X}_q \times \mathcal{C}(\overline{f})$. Then, using \eqref{derivada_Frechet_J} and \eqref{derivada_Frechet_G} in \eqref{expressao_multiplicador_de_Lagrange_aplicada} we conclude that there exists a Lagrange multiplier $\xi = (\lambda,\eta) \in L^{q'}(Q) \times L^{q'}(Q)$ such that, for all $(U,V,F) \in \widehat{X}_q \times \widehat{X}_q \times \mathcal{C}(\overline{f})$, we have
    \begin{equation} \label{desigualdade_multiplicador_de_Lagrange}
      \begin{array}{l}
         \D{\int_0^T \int_{\Omega}} (g_\lambda U + g_\eta V)\, dx \, dt 
         + \gamma_f \D{\int_0^T \int_{\Omega_c}} sgn(\overline{f}) \norm{\overline{f}}{}^{q-1} F \, dx \, dt 
          \\[6pt]
        - \D{\int_0^T \int_{\Omega}} \Big ( \partial_t U - \Delta U + \nabla \cdot (U \nabla \overline{v}) 
         + \nabla \cdot (\overline{u} \nabla V) \Big ) \lambda \, dx \, dt
         \\[6pt]
        - \D{\int_0^T \int_{\Omega}} \Big ( \partial_t V - \Delta V + s \overline{u}^{s-1} U \overline{v} + \overline{u}^s V 
        - \overline{f} V \ 1_{\Omega_c} - F \overline{v} \ 1_{\Omega_c} \Big ) \eta \, dx \, dt \geq 0
      \end{array}
    \end{equation}
    Since \eqref{desigualdade_multiplicador_de_Lagrange} is valid for all $(U,V,F) \in \widehat{X}_q \times \widehat{X}_q \times \mathcal{C}(\overline{f})$, we can deduce the optimality system \eqref{problema_adjunto_ao_linearizado} and the optimality condition \eqref{condicao_otimalidade_fraco_F}. In fact, since $\widehat{X}_q$ is a vectorial space, if $U,V \in \widehat{X}_q$ then $-U,-V \in \widehat{X}_q$. With this in mind, if we take $F = 0$ in \eqref{desigualdade_multiplicador_de_Lagrange} we obtain \eqref{adjoint-1}, \eqref{adjoint-2} and 
    \eqref{condicao_otimalidade_fraco_F}, respectively.

\
    
    \noindent{\bf Step 2: Uniqueness}

\vspace{3pt}
    
    Now, to prove the uniqueness, we suppose two possible (very weak) Lagrange multipliers
    
    \noindent
     $(\lambda_1,\eta_1), (\lambda_2,\eta_2) \in L^{q'}(Q) \times L^{q'}(Q)$ satisfying \eqref{adjoint-1} and \eqref{adjoint-2}. Let $(\widetilde{\lambda}, \widetilde{\eta}) = (\lambda_2,\eta_2) - (\lambda_1,\eta_1)$, we will prove that $\widetilde{\lambda} = \widetilde{\eta} = 0$. Subtracting the equation satisfied by $(\lambda_1,\eta_1)$ and  $(\lambda_2,\eta_2)$, then $(\widetilde{\lambda}, \widetilde{\eta})$ satisfies
    \begin{equation} \label{equacao_otimalidade_lambda_eta_barra_U}
      \begin{array}{l}
        \D{\int_0^T \int_{\Omega}} \Big ( \partial_t U - \Delta U + \nabla \cdot (U \nabla \overline{v}) \Big ) \widetilde{\lambda} + s \overline{u}^{s-1} U \overline{v} \ \widetilde{\eta} \ dx \ dt = 0, \ \forall U \in \widehat{X}_q,
      \end{array}
    \end{equation}
    \vspace{-12pt}
    \begin{equation} \label{equacao_otimalidade_lambda_eta_barra_V}
      \begin{array}{l}
        \D{\int_0^T \int_{\Omega}} \Big ( \partial_t V - \Delta V + \overline{u}^s V - \overline{f} V 1_{\Omega_c} \Big ) \widetilde{\eta} + \nabla \cdot (\overline{u} \nabla V) \widetilde{\lambda} \ dx \ dt = 0,
      \end{array}
    \end{equation}
    for all $V \in \widehat{X}_q$. Summing \eqref{equacao_otimalidade_lambda_eta_barra_U} and \eqref{equacao_otimalidade_lambda_eta_barra_V} we obtain
    \begin{equation} \label{equacao_otimalidade_lambda_eta_barra}
      \begin{array}{l}
        \D{\int_0^T \int_{\Omega}} \Big ( \partial_t U - \Delta U + \nabla \cdot (U \nabla \overline{v}) + \nabla \cdot (\overline{u} \nabla V) \Big ) \widetilde{\lambda} \ dx \ dt \\[6pt]
        + \D{\int_0^T \int_{\Omega}} \Big ( \partial_t V - \Delta V + \overline{u}^s V + s \overline{u}^{s-1} U \overline{v} - \overline{f} V 1_{\Omega_c} \Big ) \widetilde{\eta} \ dx \ dt = 0,
      \end{array}
    \end{equation}
    for all $(U,V) \in \widehat{X}_q \times \widehat{X}_q$. Now let $g_U = sgn(\tilde{\lambda}) \norm{\tilde{\lambda}}{}^{1/(q-1)}$ and $g_V = sgn(\widetilde{\eta}) \norm{\widetilde{\eta}}{}^{1/(q-1)}$. Since $(\widetilde{\lambda}, \widetilde{\eta}) \in L^{q'}(Q) \times L^{q'}(Q)$, with $q' = q/(q-1)$, we have $g_U, g_V \in L^q(Q)$. Take $(U,V) \in \widehat{X}_q \times \widehat{X}_q$ as the unique strong solution of \eqref{problema_linearizado} for this choice of $(g_U,g_V)$, therefore, from equation \eqref{equacao_otimalidade_lambda_eta_barra} we obtain $\norma{\tilde{\lambda}}{L^{q'}(Q)}^{q'} + \norma{\widetilde{\eta}}{L^{q'}(Q)}^{q'} = 0,$ which implies that $\widetilde{\lambda} = \widetilde{\eta} = 0$.
  
  \subsection{Proof of Theorem \ref{teo_regularidade_adicional_multiplicadores_de_Lagrange}}
    \label{subsec: regularidade_adicional_multiplicadores_de_Lagrange}

\
    
    \noindent {\bf Case $\boldsymbol{g_{\lambda} \in L^p(Q)}$, with $\boldsymbol{p \in [10/9,10/7)}$:} Let $\widetilde{t} = T - t$ and $\widetilde\lambda(\widetilde{t})=\lambda(t),$ $\widetilde\eta(\widetilde{t})=\eta(t)$, then problem \eqref{problema_adjunto_ao_linearizado} is equivalent to
    \begin{equation} \label{problema_adjunto_ao_linearizado_variavel_t_barra}
      \left \{
      \begin{array}{l}
        \partial_{\widetilde{t}} \widetilde\lambda - \Delta \widetilde\lambda = \nabla \overline{v} \cdot \nabla \widetilde\lambda - s \overline{u}^{s-1} \overline{v} \widetilde\eta + g_\lambda
        \\[6pt]
        \partial_{\widetilde{t}} \widetilde\eta - \Delta \widetilde\eta = - \overline{u}^s \widetilde\eta + \overline{f} \widetilde\eta \ 1_{\Omega_c} - \nabla \cdot (\overline{u} \nabla \widetilde\lambda) + g_\eta
         \\[6pt]
        \partial_{\boldsymbol{n}} \tilde\lambda |_{\Gamma} = \partial_{\boldsymbol{n}} \widetilde\eta |_{\Gamma} = 0, \ \tilde\lambda(0,x) = \widetilde\eta(0,x) = 0. 
      \end{array}
      \right.
    \end{equation}
    Then, applying Theorem \ref{teo_problema_linearizado_geral}, case $1b$, 
    with $(U,V) = (\widetilde\eta, \widetilde\lambda)$, $a_1 = \overline{u}^s - \overline{f} \ 1_{\Omega_c}$, $b_1 = 0$, $\Vec{c}_1 = 0$, $d = \overline{u}$, $g_U = g_\eta$, 
    $a_2 = 0$, $b_2 = s \overline{u}^{s-1} \overline{v}$, $\Vec{c}_2 = \nabla \overline{v}$ and 
    $g_V = g_\lambda$, 
    we conclude that there is a very weak solution $(\widetilde{\lambda}, \widetilde{\eta}) \in L^2(Q) \times L^2(Q)$ of \eqref{problema_adjunto_ao_linearizado_variavel_t_barra} and, therefore, the corresponding $({\lambda},{\eta}) \in L^2(Q) \times L^2(Q)$ is the very weak solution of \eqref{problema_adjunto_ao_linearizado} (since $q > 5/2 > 2$ we have $q' < 2$ and hence $({\lambda},{\eta}) \in L^{q'}(Q) \times L^{q'}(Q)$). Then, from the uniqueness result of Theorem \ref{teo_existencia_multiplicador_de_Lagrange_aplicado} we conclude the proof.
    
    \vspace{6pt}
    
    \noindent {\bf Case $\boldsymbol{g_{\lambda} \in L^p(Q)}$, with $\boldsymbol{p \in [10/7,2]}$:} Using the same argument of the previous case, this time applying Theorem~\ref{teo_problema_linearizado_geral}, case $1a$, with $g_U = g_\eta
    \in L^2(Q) \hookrightarrow L^{10/7}(Q)$ and $g_V = g_\lambda
     \in L^p(Q) \hookrightarrow L^{10/7}(Q)$, we conclude that the Lagrange multiplier $(\lambda,\eta)$ furnished by Theorem \ref{teo_existencia_multiplicador_de_Lagrange_aplicado} is a weak solution of \eqref{problema_adjunto_ao_linearizado} with regularity $(\lambda,\eta) \in W_2 \times W_2$. Now we enhance the regularity of $(\lambda,\eta)$ by means of a bootstrap procedure analogous to the one that was used in the proof of Theorem~\ref{teo_ponto_regular}. We first enhance the regularity of $\lambda$. Since in the right hand side of the $\lambda$-equation we have 
     $-s \overline{u}^{s-1} \overline{v} \eta + g_\lambda
      \in L^{10/3}(Q) + L^p(Q)=L^p(Q)$ (because $p \leq 2$), we apply the procedure and conclude that $\lambda \in X_p$.

    Next we apply the bootstrap to the $\eta$-equation. Since in the right hand side of the $\eta$-equation we have $- \nabla \cdot (\overline{u} \nabla \lambda) + g_\eta
    \in L^p(Q) + L^2(Q)=L^p(Q)$ ($p \leq 2$), we apply the procedure and then $\eta \in X_p$, finishing the proof.


\appendix

\section{Existence of solution for a general linear system}\label{appendixA}

  To prove the  existence of solution to the linearized problem \eqref{problema_linearizado} we introduce 
  the following general prototype of a linearized problem related to chemotaxis models:
  \begin{equation} \label{problema_linearizado_geral}
    \left \{
    \begin{array}{l}
      \partial_t U - \Delta U + a_1 U + b_1 V + \nabla \cdot (U \Vec{c}_1) + \nabla \cdot (d \nabla V) = g_U, \\
      \partial_t V - \Delta V + a_2 V + b_2 U + \Vec{c}_2 \cdot \nabla V = g_V, \\
      \partial_{\boldsymbol{n}} U |_{\Gamma} = \partial_{\boldsymbol{n}} V |_{\Gamma} = 0, \ U(0,x) = V(0,x) = 0,
    \end{array}
    \right.
  \end{equation}
  where the coefficients $a_i, b_i, d$ and $\Vec{c}_i$  are data defined in $Q$. 
  The study of \eqref{problema_linearizado_geral} will be also useful to prove  regularity of the Lagrange multiplier associated to a local optimal solution.

  In the following theorem we use the weak solutions Banach space $W_2$ defined in \eqref{weak_solutions_space}. Moreover, the concepts of weak solution, strong solution and very weak solution used in this theorem are given in Definitions \ref{defi_weak_solution}, \ref{defi_strong_solution} and \ref{defi_very_weak_solution}, respectively.
  
  \begin{theorem} \label{teo_problema_linearizado_geral}
     Let $a_i \in L^{5/2}(Q)$ and $\Vec{c}_i \in L^5(Q)^3$ with $\Vec{c}_i \cdot \Vec{n} |_{\Gamma} = 0$ for $i = 1, 2$.
    \begin{enumerate}
      \item if $b_i \in L^{5/2}(Q)$ and $d \in L^{\infty}(Q)$ we have:
      \begin{enumerate}
        \item if $g_U,g_V \in L^{10/7}(Q)$ then there is a weak solution $(U,V) \in W_2 \times W_2$ of \eqref{problema_linearizado_geral};
        \item if $g_U,g_V \in L^{10/9}(Q)$ and $\nabla d \in L^5(Q)^3$ then there is a very weak solution $(U,V) \in L^2(Q) \times L^2(Q)$ of \eqref{problema_linearizado_geral};
      \end{enumerate}
      \item if $b_1 \in L^{5/3}(Q)$ and $b_2, d \in L^5(Q)$ we have:
      \begin{enumerate}
        \item if $g_U \in L^{10/7}(Q)$ and $g_V \in L^2(Q)$ then there is a weak-strong solution $(U,V) \in W_2 \times X_2$ of \eqref{problema_linearizado_geral};
        \item if $g_U \in L^{10/9}(Q)$ and $g_V \in L^{10/7}(Q)$ then there is a very weak-weak solution $(U,V) \in L^2(Q) \times W_2$ of \eqref{problema_linearizado_geral}.
      \end{enumerate}
    \end{enumerate}
  \end{theorem}
\begin{proof}[\bf Proof]
  We will use the Galerkin method. Let $\{ \varphi_m \}$ be the basis of $H^1 $ of functions satisfying
    \begin{equation*}
      - \Delta \varphi_m + \varphi_m = \lambda_m \varphi_m, \ \partial_{\boldsymbol{n}} \varphi_m |_{\Gamma} = 0,
    \end{equation*}
    for each $m \in \mathbb{N}$, and define $X^n=\langle \varphi_1,\cdots,\varphi_n  \rangle$ 
    Let $a_i^n, b_i^n, d^n \in C^{\infty}_c(\mathbb{R} \times \mathbb{R}^3)$ and $\Vec{c}_i^{\hspace{2pt}n} \in C^{\infty}_c(\mathbb{R} \times \mathbb{R}^3)^3$ be mollifier regularizations of $a_i, b_i, d$ and $\Vec{c}_i$ such that the following strong convergences hold
    \begin{equation*}
      a_i^n \rightarrow a_i \mbox{  in } L^{5/2}(Q), \ 
       \Vec{c}_i^{\hspace{2pt}n} \rightarrow \Vec{c}_i \mbox{  in } L^5(Q)^3,\mbox{ for } i=1,2,
    \end{equation*}
    Moreover, in the case of item $(1)$ we have 
    \begin{equation*}
      b_i^n \rightarrow b_i \mbox{ strongly in } L^{5/2}(Q), \mbox{ for } i=1,2,
    \end{equation*}
    \begin{equation*}
      d^n \mbox{ is bounded in } L^{\infty}(Q) \mbox{ and converges to } d \mbox{ strongly in } L^p(Q), \mbox{ for any } p \in [1,\infty),
    \end{equation*}
    with
      $d^n \rightarrow d \mbox{ strongly in } L^5(W^{1,5})$
    in the case of item $(1b)$. On the other hand, in the case of item $(2)$, we have the strong convergences
    $
      b_1^n \rightarrow b_1 \mbox{  in } L^{5/3}(Q),\ 
       b_2^n \rightarrow b_2 \mbox{  in } L^5(Q)$ and 
       $d^n \rightarrow d \mbox{  in } L^5(Q)
    $. 
    We look for Galerkin solutions $(U_n,V_n)$ of the form $U_n(t,x) = \D{\sum_{j=1}^n{g^n_j(t) \varphi_j(x)}}$ and $V_n(t,x) = \D{\sum_{j=1}^n{h^n_j(t) \varphi_j(x)}}$ such that
    \begin{align}
      & \prodl{\partial_t U_n}{\varphi} 
      + \prodl{\nabla U_n - U_n \Vec{c}_1^{\hspace{2pt}n} -d^n \nabla V_n}{\nabla \varphi}
       + \prodl{a_1^n U_n + b_1^n V_n}{\varphi} \label{sistema_aproximado_U_n} = \prodl{g_U}{\varphi}, \\
      & \prodl{\partial_t V_n}{\varphi} 
      + \prodl{\nabla V_n}{\nabla\varphi} 
      + \prodl{a_2^n V_n+b_2^n U_n+\Vec{c}_2^{\hspace{2pt}n} \cdot \nabla V}{\varphi} 
      = \prodl{g_V}{\varphi}, \label{sistema_aproximado_V_n} \\
      & U_n(0,x) = V_n(0,x) = 0, \label{cond_iniciais_aproximado}
    \end{align}
    for all $\varphi \in X^n$. From the results on linear ordinary differential systems with smooth coefficients, we have the existence and uniqueness of global classical solution $(U_n,V_n) \in C^1([0,T];X^n \times X^n)$ satisfying \eqref{sistema_aproximado_U_n}-\eqref{cond_iniciais_aproximado}, for each $n \in \mathbb{N}$. Next we obtain \emph{a priori} estimates for $(U_n,V_n)$ that we will use to pass to the limit as $n \to \infty$. Now deal with each case of the theorem.
    
\vspace{3pt}
    
    \noindent {\bf Case $(\boldsymbol{1a})$:} We begin by taking $\varphi = U_n \in X^n$ in \eqref{sistema_aproximado_U_n} and obtain
    \begin{equation} \label{estim-1}
 \begin{array}{l}
  \displaystyle     \frac12 \frac{d}{dt} \norma{U_n}{L^2}^2 + \norma{\nabla U_n}{L^2}^2 \leq \norma{a_1^n}{L^{5/2}} \norma{U_n}{L^{10/3}}^2 + \norma{b_1^n}{L^{5/2}} \norma{V_n}{L^{10/3}} \norma{U_n}{L^{10/3}}
       \\[6pt]
  \D   \quad  + \norma{U_n}{L^{10/3}} \norma{\Vec{c}_1^{\hspace{2pt}n}}{L^5} \norma{\nabla U_n}{L^2} + \norma{d^n}{L^{\infty}} \norma{\nabla V_n}{L^2} \norma{\nabla U_n}{L^2} + \norma{g_U}{L^{10/7}} \norma{U_n}{L^{10/3}}
    \end{array}
\end{equation}
    Next, using \eqref{interpolacao_norma_10/3} and Young's inequality, we bound the last term as
    \begin{align*}
      & \norma{g_U}{L^{10/7}} \norma{U_n}{L^{10/3}} 
       \leq C_1 \norma{g_U}{L^{10/7}}^{2/7} \norma{U_n}{L^2}^{2/5} \norma{g_U}{L^{10/7}}^{5/7} \norma{U_n}{H^1}^{3/5} \\
      & \leq C_2 \norma{g_U}{L^{10/7}}^{10/7} \norma{U_n}{L^2}^2 + C_4 \norma{g_U}{L^{10/7}}^{10/7} + C_5 \norma{U_n}{H^1}^2.
    \end{align*}
    Then, applying the properties of the mollified sequences, the interpolation inequality \eqref{interpolacao_norma_10/3} and Young's inequality with the appropriate weights, we conclude that 
    there are $C, \tilde{\beta} > 0$ such that
    \begin{equation} \label{estimativa_a_priori_U_n_weak-weak}
      \begin{array}{l}
        \dfrac{1}{2} \dfrac{d}{dt} \norma{U_n}{L^2}^2 + \tilde{\beta} \norma{\nabla U_n}{L^2}^2 \leq C (\norma{a_1}{L^{5/2}}^{5/2} + \norma{b_1}{L^{5/2}}^{5/2} + \norma{\Vec{c}_1}{L^5}^5 + \norma{g_U}{L^{10/7}}^{10/7}) \norma{U_n}{L^2}^2 
         \\[6pt]
        + C \norma{V_n}{L^2}^2 + C \norma{g_U}{L^{10/7}}^{10/7} + C(\norma{d}{L^{\infty}}^2 + 1) \norma{\nabla V_n}{L^2}^2.
      \end{array}
    \end{equation}
    Now we take $\varphi = V_n \in X^n$ in \eqref{sistema_aproximado_V_n}, which gives us
    \begin{align*}
      & \frac{1}{2} \frac{d}{dt} \norma{V_n}{L^2}^2 + \norma{\nabla V_n}{L^2}^2 \leq \norma{a_2^n}{L^{5/2}} \norma{V_n}{L^{10/3}}^2 
      + \norma{b_2^n}{L^{5/2}} \norma{U_n}{L^{10/3}} \norma{V_n}{L^{10/3}} \\
      & + \norma{\Vec{c}_2^{\hspace{2pt}n}}{L^5} \norma{\nabla V_n}{L^2} \norma{V_n}{L^{10/3}} + \norma{g_V}{L^{10/7}} \norma{V_n}{L^{10/3}}.
    \end{align*}
    Applying the properties of the mollified sequences, the interpolation inequalities \eqref{interpolacao_norma_10/3} and Young's inequality with the appropriate weights, we conclude that for any $\delta > 0$ there are $C, \tilde{\beta} > 0$ (dependent on $\delta$) such that
    \begin{equation} \label{estimativa_a_priori_V_n_weak-weak}
      \begin{array}{l}
        \dfrac{1}{2} \dfrac{d}{dt} \norma{V_n}{L^2}^2 + \tilde{\beta} \norma{\nabla V_n}{L^2}^2 \leq C \norma{U_n}{L^2}^2 
          + C \norma{g_V}{L^{10/7}}^{10/7} + \delta \norma{\nabla U_n}{L^2}^2
          \\[6pt]
        + C (\norma{a_2}{L^{5/2}}^{5/2} + \norma{b_2}{L^{5/2}}^{5/2} + \norma{\Vec{c}_2}{L^5}^5 + \norma{g_V}{L^{10/7}}^{10/7}) \norma{V_n}{L^2}^2 
      \end{array}
    \end{equation}
    Summing \eqref{estimativa_a_priori_U_n_weak-weak} and $C_0$ times \eqref{estimativa_a_priori_V_n_weak-weak} and choosing $C_0 = 2 C(\norma{d}{L^{\infty}}^2 + 1)/\tilde{\beta}$ and $\delta > 0$ small enough, then the terms $\delta \norma{\nabla U_n}{L^2}^2$ and $C(\norma{d}{L^{\infty}}^2 + 1) \norma{\nabla V_n}{L^2}^2$ can be absorbed and we conclude that there is $\beta > 0$ such that
    \begin{equation} \label{estimativa_a_priori_U_n_V_n_weak-weak}
      \begin{array}{l}
        \dfrac{1}{2} \dfrac{d}{dt} (\norma{U_n}{L^2}^2 + C_0 \norma{V_n}{L^2}^2) + \beta (\norma{\nabla U_n}{L^2}^2 + \norma{\nabla V_n}{L^2}^2) \\[6pt]
        \leq C (\norma{a_1}{L^{5/2}}^{5/2} + \norma{b_1}{L^{5/2}}^{5/2} + \norma{\Vec{c}_1}{L^5}^5   + \norma{g_U}{L^{10/7}}^{10/7} + 1) \norma{U_n}{L^2}^2
        \\[6pt]
       + C (\norma{a_2}{L^{5/2}}^{5/2} + \norma{b_2}{L^{5/2}}^{5/2} 
        + \norma{\Vec{c}_2}{L^5}^5 + \norma{g_V}{L^{10/7}}^{10/7} + 1) \norma{V_n}{L^2}^2 + C \norma{g_U}{L^{10/7}}^{10/7} + C \norma{g_V}{L^{10/7}}^{10/7}.
      \end{array}
    \end{equation}
    Since $\norma{a_1}{L^{5/2}}^{5/2}$, $\norma{b_1}{L^{5/2}}^{5/2}$, $\norma{\Vec{c}_1}{L^5}^5$, $\norma{g_U}{L^{10/7}}^{10/7}$, $\norma{a_2}{L^{5/2}}^{5/2}$, $\norma{b_2}{L^5}^5$, $\norma{\Vec{c}_2}{L^5}^5$, $\norma{g_V}{L^{10/7}}^{10/7} \in L^1(0,T)$, we are able to apply Gronwall's Lemma to \eqref{estimativa_a_priori_U_n_V_n_weak-weak} and conclude that
    \begin{equation*}
      (U_n, V_n) \mbox{ is bounded in } ( L^{\infty}(L^2) \cap L^2(H^1) )^2.
    \end{equation*}
    Using this bound in the equations \eqref{sistema_aproximado_U_n} and \eqref{sistema_aproximado_V_n} we also obtain $n$-independent bounds for $\partial_t U_n$ and $\partial_t V_n$, which leads us to
    \begin{equation} \label{limitacao_W_2_W_2}
      (U_n, V_n) \mbox{ is bounded in } W_2 \times W_2.
    \end{equation}
    Next, we skip the standard procedures of the application of the Galerkin's method to linear equations and state that with \eqref{limitacao_W_2_W_2} we are able to pass to the limit as $n \to \infty$ in \eqref{sistema_aproximado_U_n} and \eqref{sistema_aproximado_V_n}, concluding that there is $(U,V) \in W_2 \times W_2$ solution of problem \eqref{problema_linearizado_geral}.
    
    \vspace{3pt}
    
    \noindent {\bf Case $(\boldsymbol{1b})$:} The $n$-independent \emph{a priori} estimates for this case are similar to those of the case $(1a)$, but now,
    we take as test functions $\varphi = (-\Delta + I)^{-1} U_n \in X^n$ in \eqref{sistema_aproximado_U_n} and $\varphi = (-\Delta + I)^{-1} V_n \in X^n$ in \eqref{sistema_aproximado_V_n}, where $\Phi = (-\Delta + I)^{-1} \phi$ is  defined as the solution of 
    \begin{equation*}
      - \Delta \Phi + \Phi = \phi, \ \partial_{\boldsymbol{n}} \Phi |_{\Gamma} = 0.
    \end{equation*}
    We also use the fact that there is a constant $C > 0$ such that $\norma{\Phi}{H^2 } \leq C \norma{ \phi}{L^2 }, \ \forall \, \phi \in L^2.$
    
    Another relevant change is that we integrate by parts to reduce the order of the space derivatives of $U_n$ and $V_n$ in \eqref{sistema_aproximado_U_n} and \eqref{sistema_aproximado_V_n} and we highlight the term $\prodl{d^n \nabla V_n}{\nabla \varphi}$ of the \eqref{sistema_aproximado_U_n} that, in this very weak solution setting, is written as $- \prodl{V_n \nabla d^n}{\nabla \varphi} - \prodl{V_n d^n}{\Delta \varphi}$.
    
    \vspace{3pt}
    
    \noindent {\bf Case $(\boldsymbol{2a})$:} We take $\varphi = V_n - \Delta V_n \in X^n$ in \eqref{sistema_aproximado_V_n}, which gives us
    \begin{align*}
      & \frac{1}{2} \frac{d}{dt} \norma{V_n}{H^1 }^2 + \norma{\nabla V_n}{L^2 }^2 + \norma{\Delta V_n}{L^2 }^2 \leq \norma{a_2^n}{L^{5/2} } \norma{V_n}{L^{10/3} }^2 + \norma{a_2^n}{L^{5/2} } \norma{V_n}{L^{10} } \norma{\Delta V_n}{L^2 } \\
      & + \norma{b_2^n}{L^5 } \norma{U_n}{L^{10/3} } \norma{V_n}{L^2 } + \norma{b_2^n}{L^5 } \norma{U_n}{L^{10/3} } \norma{\Delta V_n}{L^2 } + C \norma{\Vec{c}_2^{\hspace{2pt}n}}{L^5 } \norma{\nabla V_n}{L^{10/3} } \norma{V_n}{L^2 } \\
      & + C \norma{\Vec{c}_2^{\hspace{2pt}n}}{L^5 } \norma{\nabla V_n}{L^{10/3} } \norma{\Delta V_n}{L^2 } + \norma{g_V}{L^2 } (\norma{V_n}{L^2 } + \norma{\Delta V_n}{L^2 }).
    \end{align*}
    Applying the properties of the mollified sequences, the interpolation inequalities \eqref{interpolacao_norma_10/3} and \eqref{interpolacao_norma_10} and Young's inequality with the appropriate weights, we conclude that for any $\delta > 0$ there are $C, \tilde{\beta} > 0$ such that
    \begin{equation} \label{estimativa_a_priori_V_n_weak-strong}
      \begin{array}{l}
        \dfrac{1}{2} \dfrac{d}{dt} \norma{V_n}{H^1 }^2 + \tilde{\beta} (\norma{\nabla V_n}{L^2 }^2 + \norma{\Delta V_n}{L^2 }^2) \leq
         C \norma{g_V}{L^2 }^2 + \delta \norma{\nabla U_n}{L^2 }^2.
         \\[6pt]
     +  C (\norma{a_2}{L^{5/2} }^{5/2}  + \norma{b_2}{L^5 }^5 + \norma{\Vec{c}_2}{L^5 }^5 + \norma{g_V}{L^2 }^2 + 1) \norma{V_n}{H^1 }^2 + C \norma{U_n}{L^2 }^2 .
      \end{array}
    \end{equation}
    Now we take $\varphi = U_n \in X^n$ in \eqref{sistema_aproximado_U_n} and obtain
    \begin{align*}
      & \frac{1}{2} \frac{d}{dt} \norma{U_n}{L^2 }^2 + \norma{\nabla U_n}{L^2 }^2 \leq \norma{a_1^n}{L^{5/2} } \norma{U_n}{L^{10/3} }^2 + \norma{b_1^n}{L^{5/3} } \norma{V_n}{L^{10} } \norma{U_n}{L^{10/3} } \\
      & + \norma{U_n}{L^{10/3} } \norma{\Vec{c}_1^{\hspace{2pt}n}}{L^5 } \norma{\nabla U_n}{L^2 } + \norma{d^n}{L^5 } \norma{\nabla V_n}{L^{10/3} } \norma{\nabla U_n}{L^2 } + \norma{g_U}{L^{10/7} } \norma{U_n}{L^{10/3} }.
    \end{align*}
    Applying the properties of the mollified sequences, the interpolation inequality \eqref{interpolacao_norma_10/3} and Young's inequality with the appropriate weights, we conclude  that for any $\delta > 0$ there are $C, \tilde{\beta} > 0$ such that
    \begin{equation} \label{estimativa_a_priori_U_n_weak-strong}
      \begin{array}{l}
        \dfrac{1}{2} \dfrac{d}{dt} \norma{U_n}{L^2 }^2 
        + \tilde{\beta} \norma{\nabla U_n}{L^2 }^2 \leq 
           C \norma{V_n}{H^1 }^2 + C \norma{g_U}{L^{10/7} }^{10/7} + \delta \norma{\Delta V_n}{L^2 }^2
         \\[6pt]
      + C (\norma{a_1}{L^{5/2} }^{5/2} + \norma{b_1}{L^{5/3} }^{5/3} + \norma{\Vec{c}_1}{L^5 }^5 + \norma{d}{L^5 }^5 + \norma{g_U}{L^{10/7} }^{10/7} + 1) \norma{U_n}{L^2 }^2 .
      \end{array}
    \end{equation}
    Summing \eqref{estimativa_a_priori_V_n_weak-strong} and \eqref{estimativa_a_priori_U_n_weak-strong} and choosing $\delta > 0$ small enough so that the terms $\delta \norma{\nabla U_n}{L^2 }^2 + \delta \norma{\Delta V_n}{L^2 }^2$ on the right hand side can be absorbed by the corresponding terms on the left hand side, we conclude that there is $\beta > 0$ such that
    \begin{equation} \label{estimativa_a_priori_U_n_V_n_weak-strong}
      \begin{array}{l}
        \dfrac{1}{2} \dfrac{d}{dt} (\norma{U_n}{L^2 }^2 + \norma{V_n}{H^1 }^2) + \beta (\norma{\nabla U_n}{L^2 }^2 + \norma{\nabla V_n}{L^2 }^2) 
         + \beta \norma{\Delta V_n}{L^2 }^2
         \\[6pt]
        \leq C (\norma{a_1}{L^{5/2} }^{5/2} + \norma{b_1}{L^{5/3} }^{5/3} + \norma{\Vec{c}_1}{L^5 }^5 
        + \norma{d}{L^5 }^5 + \norma{g_U}{L^{10/7} }^{10/7} + 1) \norma{U_n}{L^2 }^2 \\[6pt]
        + C (\norma{a_2}{L^{5/2} }^{5/2} + \norma{b_2}{L^5 }^5 + \norma{\Vec{c}_2}{L^5}^5 + \norma{g_V}{L^2 }^2 + 1) \norma{V_n}{H^1 }^2
        + C \norma{g_U}{L^{10/7} }^{10/7} + C \norma{g_V}{L^2 }^2.
      \end{array}
    \end{equation}
    Since $\norma{a_1}{L^{5/2} }^{5/2}$, $\norma{b_1}{L^{5/3} }^{5/3}$, $\norma{\Vec{c}_1}{L^5 }^5$, $\norma{d}{L^5 }^5$, $\norma{g_U}{L^{10/7} }^{10/7}$, $\norma{a_2}{L^{5/2} }^{5/2}$, $\norma{b_2}{L^5 }^5$, $\norma{\Vec{c}_2}{L^5 }^5$, $\norma{g_V}{L^2 }^2 \in L^1(0,T)$,  applying Gronwall's Lemma to \eqref{estimativa_a_priori_U_n_V_n_weak-strong},
    \begin{equation*}
      (U_n, V_n) \mbox{ is bounded in } L^{\infty}(L^2) \times L^{\infty}(H^1) \cap L^2(H^1) \times L^2(H^2).
    \end{equation*}
    Using this bound in the equations \eqref{sistema_aproximado_U_n} and \eqref{sistema_aproximado_V_n} we also obtain $n$-independent bounds for $\partial_t U_n$ and $\partial_t V_n$, which leads us to
    \begin{equation} \label{limitacao_W_2_X_2}
      (U_n, V_n) \mbox{ is bounded in } W_2 \times X_2.
    \end{equation}
    Again, we skip the standard procedures of the application of the Galerkin's method to linear equations and state that with \eqref{limitacao_W_2_X_2} we are able to pass to the limit as $n \to \infty$ in \eqref{sistema_aproximado_U_n} and \eqref{sistema_aproximado_V_n}, concluding that there is $(U,V) \in W_2 \times X_2$ solution of problem \eqref{problema_linearizado_geral}.

\vspace{3pt}
    
    \noindent {\bf Case $(\boldsymbol{2b})$:} The $n$-independent \emph{a priori} estimates for this case are similar to those of the case $2a$ and one can obtain them based on the previous cases.
\end{proof}




\bibliographystyle{amsplain}

\end{document}